\documentstyle[12pt]{article}
\newcommand{\1}{{\rm 1\!\!1}}
\newcommand{\dl}{\langle \! \langle}
\newcommand{\dr}{\rangle \! \rangle}
\newcommand{\bl}{[\![}
\newcommand{\br}{]\!]}
\newcommand{\tr}{{\rm I\!P}_{0}}
\newcommand{\trr}{{\rm I\!P}_{1}}

\newcounter{mysection}
\def\myownsection{\refstepcounter{mysection} \setcounter{equation}{0}}

\begin{document}
$\;$\\[20pt]
\begin{center}
{\bf FILTERED STOCHASTIC CALCULUS}\\[60pt]
{\sc Romuald Lenczewski}\\[40pt]
Institute of Mathematics\\ 
Wroc{\l}aw University of Technology\\
Wybrze$\dot{{\rm z}}$e Wyspia{\'n}skiego 27\\
50-370 Wroc{\l}aw, Poland\\
e-mail lenczew@im.pwr.wroc.pl\\[40pt]
\end{center}
\begin{abstract}
By introducing a {\it color filtration} to the multiplicity space
${\cal G}$, we extend the quantum It$\hat{{\rm o}}$ calculus 
on multiple symmetric Fock space $\Gamma(L^{2}({\bf R}^{+}, {\cal G}))$
to the framework of {\it filtered adapted biprocesses}.
In this new notion of adaptedness,
``classical'' time filtration makes the integrands
similar to adapted processes,
whereas ``quantum'' color filtration produces their deviations
from adaptedness.
An important feature of this calculus, which we call 
{\it filtered stochastic calculus},
is that it provides an explicit interpolation 
between the main types of calculi, regardless of the type of independence,
including freeness, Boolean independence (more generally,
$m$-freeness) as well as tensor independence.
Moreover, it shows how boson calculus is ``deformed''
by other noncommutative notions of independence.
The corresponding {\it filtered It$\hat{o}$ formula} is derived. 
Existence and uniqueness of solutions of a class of
stochastic differential equations are 
established and unitarity conditions are derived.\\[10pt]
\end{abstract}
\myownsection
\begin{center}
{\sc 1. Introduction}
\end{center}
In this paper we develop a filtered version of the quantum It$\hat{{\rm o}}$
calculus on multiple symmetric Fock spaces.
It is an extension of the Hudson-Parthasarathy calculus [H-P1] 
and its multivariate version developed by Mohari and Sinha [Mo-Si]
(see also [P]).
Apart from boson calculus, it includes many other calculi,
in particular a new version of free calculus, which
was originally developed by K\"{u}mmerer and Speicher [K-Sp] for the 
Cuntz algebra, as well as
new examples of $m$-free calculi for the 
$m$-free Brownian motions introduced in [F-L], where $m\in {\bf N}$
(for an inclusion of the calculus on the finite difference algebra
[B] see [P-Si]).

In [L2] we introduced {\it filtered random variables},
from which other random variables can be obtained 
by addition or strong limits,
regardless of the notion of independence.
In particular, this includes the three main types in the axiomatic
approach to independence ([Sp1],[S2]), 
corresponding to tensor, free [V] and Boolean products of states.
The same is true for $m$-free random variables
for all $1\leq m \leq \infty$
obtained from the {\it hierarchy of freeness} 
construction [L1] (see also [F-L] for limit theorems and
[F-L-S] for the GNS construction).

By studying the asymptotic joint distributions of their normalized 
sums in limit theorems [L2], we were led to 
{\it filtered creation, annihilation, number} and {\it time
processes} (see (2.4)-(2.7)). 
They live in a multiple symmetric Fock space $\Gamma({\cal H})$, where
${\cal H}=L^{2}({\bf R}^{+}, {\cal G})$
and ${\cal G}$ is a separable Hilbert space with a countable
fixed orthonormal basis $(e_{n})_{n\in {\bf N}}$ and
are obtained from the CCR processes by multiplying them
by canonical projections
$$
P^{(V)}:\;\;\;\Gamma({\cal H})\rightarrow 
\Gamma(L^{2}({\bf R}^{+}, {\cal G}^{(V)})),\;\;
{\cal G}^{(V)}=\bigoplus_{m\in V}{\bf C}e_{m}
$$
where $V\subseteq {\bf N}$ and we set 
${\cal G}^{(\emptyset)}=\{0\}$. In other words, $P^{(V)}$
is the projection onto the subspace built from the vaccum vector
$\Omega$ and those copies
(or {\it colors}) of $L^{2}({\bf R}^{+})$ which are associated with
the set $V$.

This leads to 
the {\it filtered stochastic calculus} developed in this paper, where we
deal with integrals of type
$$
I^{\eta}(t)= \int_{0}^{t}FdA^{\eta}G,
$$
defined on the exponential domain,
with the integrator $A^{\eta}=(A^{\eta}_{t})_{t\geq 0}$ 
being one of the 
CCR basic integrators: $A^{(k)*}$ (creation), $A^{(k)}$ (annihilation),
$A^{(k)\circ}$ (number), associated with color $k$, or $A^{(0)}$ (time).
The integrands are 
biprocesses $F\otimes G =(F(t)\otimes G(t))_{t\geq 0}$ which
are not adapted, namely
$$
F(t)=\widetilde{F}(t)\otimes P^{(D)}, \;\; 
G(t)=\widetilde{G}(t)\otimes P^{(E)}
$$
for all $t\geq 0$, where
$D,E \subseteq {\bf N}$, according to the past-future decomposition
$\Gamma({\cal H})=\Gamma({\cal H}_{t]})\otimes \Gamma({\cal H}_{[t})$.
In other words, the identity corresponding to
``the future'' is replaced by color projections
with filters $D,E$ showing which colors are 
filtered through.
We will say that $F\otimes G$ is $(D,E)$-- {\it adapted}, whereas linear 
combinations of such biprocesses, corresponding to different filters,
will be called {\it filtered adapted}.

We arrive at the {\it filtered It$\hat{o}$ formula}, which takes
a particularly nice form. 
Namely, let $A^{\eta_{1}}$ and $A^{\eta_{2}}$ 
be CCR integrators associated with colors 
$k_{1}$ and $k_{2}$, respectively,
and let $dA^{\eta_{1}}dA^{\eta_{2}}=d\bl A^{\eta_{1}}, A^{\eta_{2}} \br$
be the result of the It$\hat{{\rm o}}$ multiplication
of the CCR differentials. Then 
all nontrivial It$\hat{{\rm o}}$ corrections for the differentials
$$
dI^{\eta_{1}}_{1}
=
G_{1}dA^{\eta_{1}}F_{1},\;\;\;
dI^{\eta_{2}}_{2}
=
F_{2}dA^{\eta_{2}}G_{2}
$$
where $G_{1}\otimes F_{1}$ and 
$F_{2} \otimes G_{2}$ are 
$(E_{1},D_{1})$-- and $(D_{2},E_{2})$-- adapted
locally square integrable biprocesses, respectively,
can be written as 
$$
dI_{1}^{\eta_{1}}dI_{2}^{\eta_{2}}
=\1_{D_{1}\cap D_{2}}(k_{1})
G_{1} d\bl A^{\eta_{1}}, A^{\eta_{2}}\br 
F_{1}F_{2}G_{2}
$$
where $\1_{A}$
is the indicator function of the set $A$. We proceed further and
study existence, uniqueness and unitarity of solutions of stochastic
differential equations.

It can be seen that the key role in this approach is 
played by the new notion of adaptedness, which exhibits
an interplay between ``classical'' {\it time filtration} and ``quantum''
{\it color filtration}.
The corresponding calculus is 
a non-trivial, but quite natural and general extension 
of boson calculus.
It includes many examples of
quantum stochastic calculi, gives new ones, 
like $m$-free calculi for all natural $m$, and shows
connections between them.

Let us recall that the construction of the hierarchy of freeness 
[L1] showed how to approximate the free product of states 
using tensor independence. 
Filtered stochastic calculus preserves
this hierarchy and shows that the $m$-free calculus exhibits
the $m$-th ``level of adaptedness''.
More importantly, it allows us to compare the $m$-free calculi
against other calculi, in particular, boson calculus, 
by measuring their ``deviations from adaptedness''. 
In order to do that, 
it is enough to give the collections of filters
associated with the calculi. 
This gives
\begin{center}
$1$-free calculus -- ${\cal P}^{(1)}=\{\emptyset , \{1\}\}$\\
\ldots\\
$m$-free calculus -- ${\cal P}^{(m)}=\{
\emptyset , \{1\}, \ldots  \{1, \ldots , m\}\}$\\
\ldots\\
free calculus -- ${\cal P}^{(\infty)}=
\{\emptyset , \{1\}, \ldots  \{1, \ldots , m\}, \ldots \}$
\end{center}
for the ``minimal'' formulation, i.e. when the integrands
belong to the *-algebra generated by the corresponding fundamental
processes (if we add the unit or study unitarity, 
we need to add the filter ${\bf N}$). 

These can be compared against the two extreme cases
\begin{center}
$\Omega$-adapted calculus -- ${\cal P}^{(0)}=\{\emptyset\}$\\[5pt]
boson calculus -- ${\cal P}=\{{\bf N}\}$,
\end{center}
i.e. the ``least adapted'' calculus studied in [Vi] and [Be] 
associated with 
the projection $P_{\Omega}$ on the zero-particle space, and
the boson calculus -- 
the ``most adapted'' calculus associated with 
the projection $P^{({\bf N})}=I$ on the whole space.
The ``least adapted'' calculus from the hierarchy of $m$-free
calculi is the Boolean (or $1$-free) calculus.
In turn, the $m$-free calculus corresponds to
a mixture of $m+1$ types of adaptedness, 
whereas the free calculus -- to a mixture of infinitely many
types of adaptedness. 
Of course, this picture holds on $\Gamma({\cal H})$. When we restrict
ourselves to suitable proper subspaces of $\Gamma({\cal H})$, for instance,
to $m$-free Fock spaces, one can construct calculi which
become adapted on those subspaces.

Let us mention here other unified approaches to quantum 
stochastic calculus. 
A representation-free stochastic calculus was 
presented by Accardi, Fagnola and Quaegebeur in 
[Ac-Fa-Qu] and [Fa].
We would like to mention here that it seems possible to treat the filtered
calculus in a similar manner by proving semimartingale inequalities
for $(D,E)$-adapted bi-processes.
Non-causal approaches
to stochastic calculus were developed by Lindsay [Li] and Belavkin [Bel]. 
For other calculi, see [Ba-St-Wi], [Ap-H], [H-P2],
[Bi-Sp], [Me], [At-Li], [Ma], [S1], [Sp2].
Our approach to calculus is closest to that of
Parthasarathy and Sinha who realized [P-Si]
free fundamental processes as stochastic
integrals of non-adapted processes in boson calculus.
However, we can go much further and treat in a unified manner 
the It$\hat{{\rm o}}$ formula as well as 
stochastic differential equations and
unitary evolutions. 

We would like to stress that our approach
is not restricted to the stochastic calculus. 
On the contrary, it provides a unified treatment of 
such elements of noncommutative probability as
(i) product states,
(ii) limit theorems,
(iii) quantum It$\hat{{\rm o}}$ calculus,
(iv) stochastic differential equations,
of which the first two were treated in [L2]. At the same time, 
it provides a very concrete mathematical framework and exhibits connections
with many other models, some of which were mentioned above.
Although we have taken the most traditional 
(Hudson-Parthasarathy) approach to
calculus, the core of our approach seems more
universal. The main idea boils down to introducing the second
filtration in the 
underlying Hilbert space (the infinite tensor product of
Hilbert spaces for the first quantization and 
multiple symmetric Fock space 
for the second quantization).
It seems likely that in other models
one can use the same idea,
which could be of importance in our understanding how
noncommutative notions of independence ``deform'' classical
probability. \\[10pt]
\myownsection
\begin{center}
{\sc 2. Definitions and notation}
\end{center}
\begin{center}
{\it Exponential domain}
\end{center}
Let ${\cal G}$ be a separable Hilbert space with a countably infinite 
fixed orthonormal basis
$(e_{n})_{n\in {\bf N}}$. It is sometimes called the {\it multiplicity space}.
By a {\it multiple symmetric Fock space } over
${\cal K}$ we understand the symmetric Fock space over
${\cal H}=L^{2}({\bf R}^{+}, 
{\cal G})\cong L^{2}({\bf R}^{+})
\otimes {\cal G}\equiv {\cal K}\otimes {\cal G}$, namely
$$
\Gamma({\cal H})={\bf C}\Omega \oplus \bigoplus_{n=1}^{\infty}{\cal H}^{\circ n}
$$
where ${\cal H}^{\circ n}$ denotes the $n$-th symmetric tensor power
of ${\cal H}$ and $\Omega$ is the vacuum vector, 
with the scalar product given by $\langle\Omega , \Omega \rangle=1$,
$\langle \Omega , u \rangle =0$ and
$$
\langle u_{1} \circ \ldots \circ u_{n},
v_{1} \circ \ldots \circ v_{m} \rangle
=\delta_{n,m}\frac{1}{n!}
\sum_{\sigma \in {\bf S}_{n}}
\langle u_{1} , v_{\sigma(1)} \rangle
\ldots 
\langle u_{n} , v_{\sigma(n)} \rangle
$$
where
$$
u_{1} \circ \ldots \circ u_{n}=
\frac{1}{n!}\sum_{\sigma \in {\bf S}_{n}}
u_{\sigma(1)}\otimes \ldots \otimes u_{\sigma(n)}
$$
and ${\bf S}_{n}$ denotes the symmetric group of order $n$.

The exponential vectors are given by
$$
\varepsilon (u)=\bigoplus_{n=0}^{\infty}\frac{1}{\sqrt{n!}}u^{\otimes n}
$$
where $u^{\otimes 0}= \Omega$, and $u\in {\cal H}$. 
Thus, in particular,
$\varepsilon (0)=\Omega$. The linear space
${\cal E}$ spanned by exponential vectors is usually called
the exponential domain.
It is well-known that ${\cal E}$ is dense in 
the symmetric Fock space $\Gamma({\cal H})$. 
The scalar product of two exponential vectors
is given by 
$$
\langle \varepsilon (u) , \varepsilon (v) \rangle = e^{\langle u , v
\rangle}.
$$
where $u,v\in {\cal H}$.

We will use the functorial property of $\Gamma({\cal H})$
for the {\it time filtration} of the Hilbert space ${\cal H}$. Namely,
for the direct sum decomposition
$$
{\cal H}={\cal H}_{s]}\oplus {\cal H}_{[s,t]}\oplus {\cal H}_{[t},
$$
where ${\cal H}_{s]}=L^{2}([0,s])\otimes {\cal G}$,
${\cal H}_{[s,t]}=L^{2}([s,t])\otimes {\cal G}$ and
${\cal H}_{[t}= L^{2}([t, \infty))\otimes {\cal G}$, we have
$$
\Gamma=\Gamma_{s]}\otimes \Gamma_{[s,t]}\otimes \Gamma_{[t}
$$
where
$\Gamma_{s]}=\Gamma({\cal H}_{s]})$,
$\Gamma_{[s,t]}=\Gamma({\cal H}_{[s,t]})$,
$\Gamma_{[t}=  \Gamma({\cal H}_{[t})$
for any $0<s<t<\infty$.

Another direct sum decomposition of ${\cal H}$ is associated with the discrete
{\it color filtration} in ${\cal H}$. 
More generally, for arbitrary $V \in {\cal P}({\bf N})$, let
$
\Pi^{(V)}:\; {\cal H} \rightarrow {\cal H}^{(V)}
$
be the canonical projection onto
\begin{equation}
\label{2.1}
{\cal H}^{(V)}=\bigoplus_{k\in V}{\cal K}\otimes e_{k}
\end{equation}
with ${\cal H}^{(\emptyset )}=\{0\}$. We will say that $\Pi^{(V)}$
is the projection onto the subspace spanned by 
vectors of colors which are in $V$.
If $V=\{1, \ldots , r-1\}$, then a short-hand notation will be used, namely
$$
{\cal H}^{(r)}:= {\cal H}^{(\{1, \ldots , r-1\})}.
$$
The vector subspace of ${\cal H}$ spanned by all
vectors $u$ of {\it finite color support}, i.e.
$u\in {\cal H}^{(r)}$ for some $r\in {\bf N}$,
will be denoted by ${\cal H}_{0}$. 
Similarly, if $x\in \Gamma({\cal H}^{(r)})$ for some $r\in {\bf N}$,
we will also say that it is of finite color support.
Finally, let 
$P^{(V)}:\;\; \Gamma({\cal H})\rightarrow \Gamma({\cal H}^{(V)})$
denote the second quantization of $\Pi^{(V)}$, thus
$P^{(V)}\varepsilon (u)=\varepsilon (\Pi^{V}u)$.

The space $\Gamma({\cal H})$ will be extended to 
$$
\widetilde{\Gamma}({\cal H})=h_{0}\otimes \Gamma ({\cal H}),
$$
where
$h_{0}$ is a separable Hilbert space called the initial space.
The ampliations
$1\otimes A^{\eta}$ of fundamental operators $A^{\eta}$
will also be denoted by $A^{\eta}$. Moreover,
by $P^{(V)}$ we will denote the ampliations $1\otimes P^{(V)}$
and also its restrictions to $\widetilde{\Gamma}_{s]}$, 
$\widetilde{\Gamma}_{[s,t]}$ or $\widetilde{\Gamma}_{[t}$.

The exponential domain ${\cal E}_{0}$, i.e. the span
of $\varepsilon (u)$, where $u$ is locally bounded as a function of time
and is of finite color support, is then replaced by
$\widetilde{\cal E}_{0}={\rm span}\{{\cal M}_{0}\}$
where
$$
{\cal M}_{0}=\{w\otimes \varepsilon (u)\equiv
w\varepsilon (u): w\in {\cal D}_{0}, u\in {\cal H}_{0}\},
$$
${\cal D}_{0}$ is a dense subset of $h_{0}$
and $u$ is locally bounded as a function of time.
Clearly, ${\cal M}_{0}$ is total and 
$\widetilde{\cal E}_{0}$ is dense in
$\widetilde{\Gamma}({\cal H})$. 
We will also write
$$
x^{(V)}=w\varepsilon (u^{(V)}).
$$
for $x\in {\cal M}_{0}$,
where $V\in {\cal P}({\bf N})$.

The notations associated with the continuous
tensor product decompositions of 
$\widetilde{\Gamma}({\cal H})$ as well as ${\cal M}_{0}$ 
will also be standard.
For instance, 
if $x=w\varepsilon (u)\in {\cal M}_{0}$, then we will write
$$
x=x_{s]}\otimes x_{[s,t]}\otimes x_{[t},
$$
where $x_{s]}=w\varepsilon (u_{s]})$,
$x_{[s,t]}=\varepsilon (u_{[s,t]})$ and $x_{[t}=\varepsilon (u_{[t})$
for $s<t$.\\
\begin{center}
{\it Filtered fundamental processes}
\end{center}
The filtered creation, annihilation and number
operators introduced in [L2] lead to the filtered
fundamental processes. Namely, they
are expressed in terms of the canonical ones 
on $\Gamma({\cal H})$ by the following formulas:
\begin{eqnarray}
\label{2.2}
A_{t}^{(k,V)*}
&=&A_{t}^{(k)*}P^{(V)}=
a^{*}(\chi_{[0,t]}\otimes e_{k})P^{(V)}\\
\label{2.3}
A^{(k,V)}_{t}
&=&
P^{(V)}A_{t}^{(k)}= P^{(V)}a(\chi_{[0,t]}\otimes e_{k})\\
\label{2.4}
A^{(k,V)\circ}_{t}
&=&
A^{(k)\circ}_{t}P^{(V\cup\{k\})}=
\lambda(I_{[0,t]}\otimes |e_{k}\rangle \langle e_{k}|)P^{(V\cup\{k\})}\\
\label{2.5}
A^{(0,V)}_{t}
&=& A^{(0)}_{t}P^{(V)}=tP^{(V)},
\end{eqnarray}
where $k\in {\bf N}$, $V\in {\cal P}({\bf N})$,
$I_{[0,t]}$ denotes the
operator of multiplication by the characteristic function 
$\chi_{[0,t]}$ on $L^{2}({\bf R}^{+})$ and
$\lambda(T)$ is the differential
second quantization of $T$. 
The 
families of processes given by
(2.2-2.5)
will be called {\it filtered creation, annihilation, 
number} and {\it time procesess}, respectively.

Clearly, if $V={\bf N}$, equations (2.2-2.5) give CCR creation, annihilation, number and
time processes, respectively. By
$$
{\cal T}=\{(k), (k)*, (k)\circ, (0)|k \in {\bf N}\}
$$
we denote the set of indices associated with these processes.
Here, we do not follow the notation of [Mo-Si] in order to
have a clear connection with free processes.
We will also find it convenient to use the ``duals'' of 
$A^{\eta}$, $\eta\in {\cal T}$. Thus 
$A^{(k)*\dagger}=A^{(k)}$, $A^{(k)\dagger}=A^{(k)\dagger}$,
$A^{(k)\circ\dagger}=A^{(k)\circ}$, $A^{(0)\dagger}=A^{(0)}$.\\
\begin{center}
{\it $m$-free fundamental processes}
\end{center}
From the extended $m$-free fundamental operators
defined in [L2] we obtain
extended $m$-free fundamental processes
\begin{eqnarray}
\label{2.6}
l^{(m)*}_{t}
&=&
\sum_{k=1}^{m}A^{(k)*}_{t}P^{[k-1]},\\
\label{2.7}
l^{(m)}_{t}
&=&
\sum_{k=1}^{m}P^{[k-1]}A^{(k)}_{t}\\
\label{2.8}
l^{(m)\circ}_{t}
&=&
\sum_{k=1}^{m}P^{[k]}A^{(k)\circ}_{t}\\
\label{2.9}
l^{(m)\cdot}_{t}
&=&
P^{(m)}A^{(0)}_{t}
\end{eqnarray}
i.e. the
{\it extended $m$-free creation, annihilation,
number} and {\it time processes}, respectively, 
where $P^{[k-1]}=P^{(k)}-P^{(k-1)}$, $P^{[0]}=P_{\Omega}$, 
and $m\in {\bf N}^{*}={\bf N}\cup \{\infty\}$.

Note that 
\begin{eqnarray*}
P^{(k)}:\;\;&&
\widetilde{\Gamma}({\cal H})\rightarrow \widetilde{\Gamma}({\cal H}^{(k)})\\
P^{[k-1]}: \;\; && \widetilde{\Gamma}({\cal H})\rightarrow 
\widetilde{\Gamma}({\cal H}^{(k)})\ominus \widetilde{\Gamma}
({\cal H}^{(k-1)}),
\end{eqnarray*}
are the orthogonal projections on the closed
subspaces of $\widetilde{\Gamma}({\cal H})$ spanned by 
elementary tensors constructed from vectors of 
colors $0, \ldots , k-1$, and those with largest color
$k-1$, respectively.

In the sequel we will make the identifications: 
$l_{t}\equiv l_{t}^{(\infty)}$, $l^{*}_{t}\equiv l^{(\infty)*}$, 
$l^{\circ}_{t}\equiv l^{(\infty)\circ}$ and 
$l^{\cdot}_{t}\equiv l^{(\infty)\cdot}$. 
We should mention that in this paper we also identify
$l^{(m)}_{t}$, $l^{(m)*}_{t}$, $l^{(m)\circ}_{t}$, $l^{(m)\cdot}_{t}$
with their ampliations to $\widetilde{\Gamma}({\cal H})=h_{0}\otimes \Gamma$.
The ``duality'' relations read:
$l^{(m)\dagger}=l^{(m)*}$,
$l^{(m)*\dagger}=l^{(m)}$, $l^{(m)\circ \dagger}=l^{(m)\circ}$
and $l^{(m)\cdot \dagger}=l^{(m)\cdot}$. 

A shorthand notation for (\ref{2.6})-(\ref{2.9}) will be used, namely
\begin{equation}
\label{2.10}
l_{t}^{\alpha}=\sum_{(\eta, V)\sim\alpha}A_{s}^{(\eta, V)}
\end{equation}
where $(\eta,V)\sim\alpha$ means that
$A^{(\eta ,V)}$ appears on the RHS of (\ref{2.6})-(\ref{2.9}) and
$$
\alpha\in {\cal F}_{m}=\{(m),(m)*,(m)\circ ,(m)\cdot\},
$$
the set of the indices of $m$-free fundamental processes, 
for $m\in {\bf N}^{*}$. 

It was shown in [L2] that the extended $m$-free fundamental processes
approximate extended free fundamental processes as $m\rightarrow \infty$,
the latter being obtained for $m=\infty$
This holds on all of $\widetilde{\Gamma}({\cal H})$ in the case
of creation, annihilation and time processes since they have unique
bounded extensions to $\widetilde{\Gamma}({\cal H})$, or on a dense domain
in the case of number processes, 
which are unbounded on $\widetilde{\Gamma}({\cal H})$.
One should also note that if $m=\infty$, then
the restrictions of formulas (\ref{2.6})-(\ref{2.9})
to $\widetilde{\Gamma}({\cal H}^{(r)})$ are always finite sums.
\begin{center}
{\it Notation $\#$}
\end{center}
We will follow [Bi-Sp] and use the symbol $\#$ 
to write stochastic integrals and their matrix elements
in such a way that $F$ and $G$ from an integrated elementary 
biprocess $F\otimes G$ are not separated by the integrators. Thus
$$
dI=F\otimes G\# dM:= FdMG
$$
will be the differential w.r.t. $dM$ and, if
$X=\sum_{i}F_{i}\otimes G_{i}$
is a biprocess integrable w.r.t. $dM$, we will write the integrals as
$$
\int_{0}^{t}X\# dM:
=
\sum_{i}\int_{0}^{t}F_{i}dMG_{i}.
$$
When calculating matrix elements of stochastic integrals,
we will also use.
In turn, when using matrix elements, we will use 
$$
(x, F\otimes G y)\# \dl I,Q \dr:=
\langle x, FQGy \rangle 
$$
where $F\otimes G$ is a suitable stochastic biproces,
$Q$ is a projection and $x, y\in \widetilde{\cal E}_{0}$,
and extend it by linearity .\\[10pt]
\myownsection
\begin{center}
{\sc 3. Filtered adapted biprocesses}
\end{center}
By ${\cal L}(h, {\cal D})$, 
where $h$ is a separable Hilbert space and ${\cal D}$ is a dense subset
of $h$, we denote the vector space of all 
linear operators $F$ on $h$ such that 
${\cal D}\subset D(F)\cap D(F^{*})$,
where $F^{*}$ is the adjoint of $F$.
In this paper we will use ${\cal L}(h_{0},{\cal D}_{0})$ 
and ${\cal L}({\Gamma}, {\cal E}_{0})$.

Let us first specify the notion of adaptedness. 
The dependence on ${\cal E}_{0}$ and ${\cal D}_{0}$
of this notion of adaptedness
is supressed in the notation.\\
\indent{\par}
{\sc Definition 3.1.}
Let $t\geq 0$ and let $D,E \in {\cal P}({\bf N})$. 
We will say that an operator
$A\otimes B$, where $A,B\in {\cal L}(h_{0},{\cal D}_{0})\otimes
{\cal L}(\Gamma, {\cal E}_{0})$,
is $(t,D,E)$ -- {\it adapted} if
\indent{\par}(1) 
$\widetilde{\cal E}_{0}\subset D(AB)$ and
$\widetilde{\cal E}_{0}\subset D(B^{*}A^{*})$ 
\indent{\par}(2) 
$A=\widetilde{A}\otimes P^{(D)}$ and
$B=\widetilde{B}\otimes P^{(E)}$ according 
to the decomposition 
$\widetilde{\Gamma}=\widetilde{\Gamma}_{t]}\otimes \widetilde{\Gamma}_{[t}$,
\indent{\par}(3) 
$A$ and $B$ leave invariant the span of vectors of finite color support, i.e.
$$
\forall r\in {\bf N}\; \exists
p,q\in {\bf N}:\;\widetilde{B}: 
(\widetilde{\cal E}^{(r)}_{0})_{t]}\rightarrow
\widetilde{\Gamma}({\cal H}^{(p)}_{t]})\;\;{\rm and} \;\;
\widetilde{A}:\widetilde{\Gamma}({\cal H}^{(p)}_{t]})\cap D(\widetilde{A})
\rightarrow \widetilde{\Gamma}({\cal H}^{(q)}_{t]}).
$$
\indent{\par}
{\sc Definition 3.2.}
Let $D,E\in {\cal P}({\bf N})$.
By an 
{\it elementary $(D,E)$-- adapted stochastic biprocess} we will understand
a family $(F(t)\otimes G(t))_{t\geq 0}$, where
\indent{\par}
(1) $F(t)\otimes G(t)$ is $(t,D,E)$ -- adapted
for all $t\geq 0$ and such that for each $r$ of Definition 3.1, the numbers $p,q$
can be chosen the same for all $t\geq 0$,
\indent{\par}
(2) the map  $t\rightarrow F(t)G(t)x$ is strongly measurable
for all  $x\in \widetilde{\cal E}_{0}$.\\
We will often denote this biprocess by
$F\otimes G$ understanding that $F=F(t)$ and $G=G(t)$.\\
\indent{\par}
{\sc Definition 3.3.}
An elementary $(D,E)$-- adapted stochastic biprocess $F\otimes G$
will be called {\it simple} 
if there exists a partition of ${\bf R}^{+}$ given by 
$0=t_{0}<t_{1}<\ldots <t_{n}< \ldots$, where $t_{n}\uparrow \infty$, 
such that 
$$
F(t)\otimes G(t)=
\sum_{k=0}^{\infty} F(t_{k})\otimes G(t_{k})\chi_{[t_{k},t_{k+1})}(t)
$$
for any $t\in  {\bf R}^{+}$.
It will be called 
{\it regular} (or, {\it continuous}) 
if the map $t\rightarrow F(t)G(t)x$
is strongly continuous for all $x\in \widetilde{\cal E}_{0}$.\\
\indent{\par}
The vector space spanned by elementary $(D,E)$-adapted
biprocesses will be denoted by 
${\cal A}(D,E)$. 
The vector subspaces of ${\cal A}(D,E)$ 
spanned by elementary simple $(D,E)$- adapted 
biprocesses, and elementary regular $(D,E)$- adapted biprocesses,
will be denoted ${\cal S}(D,E)$ and ${\cal C}(D,E)$, respectively. 
Arbitrary elements of ${\cal A}(D,E)$, 
${\cal S}(D,E)$, or ${\cal C}(D,E)$
will be called $(D,E)$- {\it adapted},
{\it simple} $(D,E)$- {\it adapted},
and {\it regular} $(D,E)$- {\it adapted biprocesses},
respectively. \\
\indent{\par}
{\sc Definition 3.4.}
Let ${\cal P}_{0}$ be a finite subset of ${\cal P}({\bf N})$.
A finite linear combination of the form
$$
X=\sum_{i}F_{i}\otimes G_{i}
$$
where $F_{i}\otimes G_{i}\in {\cal A}(D_{i},E_{i})$, $D_{i},E_{i}\in 
{\cal P}_{0}$ for all $i$, will be called a $({\cal P}_{0},{\cal P}_{0})$--
{\it adapted} stochastic biprocess. The vector space spanned by 
$({\cal P}_{0},{\cal P}_{0})$-- adapted stochastic 
biprocesses will be denoted by ${\cal A}({\cal P}_{0},{\cal P}_{0})$.
Any element of the algebraic direct sum
$$
{\cal A}=\bigoplus_{D,E\in {\cal P}({\bf N})}{\cal A}(D,E)
$$
will be called a {\it filtered adapted} stochastic biprocess.\\
\indent{\par}
A family $F=(F(t))_{t\geq 0}\equiv (F_{t})_{t\geq 0}$ of operators
from ${\cal L}(h_{0}, {\cal D}_{0})\otimes {\cal L}(\Gamma , {\cal E}_{0})$
will be called a $V$-{\it adapted stochastic process}, where 
$V\in {\cal P}({\bf N})$
if $F\otimes I$ and $I\otimes F$ 
are $(V,{\bf N})$-- and 
$({\bf N},V)$ --adapted stochastic biprocesses, respectively.
The vector space spanned by $V$-adapted
processes will be denoted by 
${\cal A}(V)$.
Note that if $F$ is ${\bf N}$ -- adapted,
then it is adapted in the usual ([H-P1]) sense. 
Note also that if $F\otimes G\in {\cal A}(D,E)$, then
$FG\in {\cal A}(D\cap E)$.
A process $F$ will be called simple,
or regular if the biprocesses $F\otimes I$, $I\otimes F$
are simple, or regular, respectively. 
A stochastic process $F$ will be called
${\cal P}_{0}$--{\it adapted} if 
$$
F(t)=\sum_{V\in {\cal P}_{0}}F_{V}(t)
$$
for all $t$, where $F_{V}\in {\cal A}(V)$ for all $V\in {\cal P}_{0}$.

Filtered fundamental processes $A^{(k,V)}$, $A^{(k,V)*}$,
$A^{(k,V)\circ}$ and $A^{(0,V)}$ are 
natural examples of 
$V$-adapted regular stochastic processes 
for all $k \in {\bf N}$, $V\in {\cal P}({\bf N})$.
In view of (2.2)=(2.5), however, 
the integrals with filtered fundamental processes
as integrators
can be expressed as stochastic integrals with boson integrators.
Namely, given a filtered adapted
biprocess $X=\sum_{i}F_{i}\otimes G_{i}$, 
$\eta \in {\cal T}$ and $V\in {\cal P}({\bf N})$,
there exists a filtered adapted biprocess $X[\eta , V]$
such that 
\begin{equation}
\label{3.1}
\int_{0}^{t}X\# dA^{(\eta,V)}
=
\int_{0}^{t}X[\eta ,V]\# dA^{\eta}
\end{equation}
on $\widetilde{\cal E}_{0}$, 
provided the integral on the RHS exists, where
\begin{equation}
\label{3.2}
X[\eta, V]=
\left\{
\begin{array}{lll}
\sum_{i}F_{i}P^{(V)}\otimes G_{i} & {\rm if} & \eta =(k)\\
\sum_{i}F_{i}\otimes P^{(V)}G_{i} & {\rm if} & \eta = (k)*\\
\sum_{i}F_{i}P^{(V\cup\{k\})}\otimes G_{i} & {\rm if} &\eta =(k)\circ\\
\sum_{i}F_{i}P^{(V)}\otimes G_{i} & {\rm if} & \eta =(0)
\end{array}
\right..
\end{equation}
\indent{\par}
Therefore, our study will concentrate on stochastic integrals 
w.r.t. the CCR processes.
Note that in the case of number and time operators, the 
projections commute with $A^{\eta}$, 
so one can also flip the projections to the other side
of the tensor product.
\\[10pt]
\newpage
\myownsection
\begin{center}
{\sc 4. Filtered fundamental lemmas}
\end{center}
Let us begin with the definition of stochastic integrals of 
elementary simple biprocesses.\\
\indent{\par}
{\sc Definition 4.1.}
Let $F\otimes G \in {\cal S}(D,E)$, where $D,E\in {\cal P}({\bf N})$,
be given by Definition 3.3. For any $\eta \in {\cal T}$, define
$$
I^{\eta}(t)=
\int_{0}^{t}
F\otimes G \# dA^{\eta}=
\sum_{k=1}^{n+1}F(t_{k-1}) \otimes G(t_{k-1}) \#
(A_{t_{k}\wedge t}^{\eta} - A_{t_{k-1}\wedge t}^{\eta})
$$
on $\widetilde{\cal E}_{0}$, where $t_{n}\leq t <t_{n+1}$. \\
\indent{\par}
Note that this definition does not depend on the partition of ${\bf R}^{+}$
in the sense that one can take a refinement of the given partition to obtain
the same result. Therefore, for convenience, we will fix $t$
and from now on assume that $t=t_{n}$.

It is convenient to introduce some notation for complex-valued measures
which appear in boson multivariate calculus. Thus, for given
$u,v\in {\cal H}_{0}$ (suppressed in the notation), let
\begin{equation}
\label{4.1}
\mu^{\eta}([s,t])=
\left\{
\begin{array}{ccl}
\int_{s}^{t}v^{(k)}(r) dr & {\rm if} & \eta=(k)\\
\int_{s}^{t}\bar{u}^{(k)}(r) dr & {\rm if} & \eta =(k)*\\
\int_{s}^{t}\bar{u}^{(k)}(r)v^{(k)}(r) dr &{\rm if} &\eta=(k)\circ\\
t-s & {\rm if} &\eta =(0)
\end{array}
\right.
\end{equation}
be the measures associated with the annihilation, creation, number
and time processes, respectively, of boson multivariate calculus.
They are all absolutely continuous with respect to the Lebesgue 
measure.\\
\indent{\par}
{\sc Lemma 4.2.}
{\it Let $0\leq s \leq t$, $x=w\varepsilon (u)$,
$y= z\varepsilon (v)$, where $w,y\in {\cal D}_{0}$,
$u,v \in {\cal H}_{0}$ and let 
$F\otimes G\in {\cal S}(D,E)$, where $D,E\in {\cal P}({\bf N})$.
Then}
$$
\langle 
x, 
I^{\eta}(t) 
y
\rangle
=\int_{0}^{t}
\langle
x,
F(s)G(s) 
y 
\rangle
d{\mu}^{\eta}_{D,E}(s)
$$
{\it where $\mu^{\eta}_{D,E}=\1^{\eta}_{D,E}\mu^{\eta}$ and}
\begin{equation}
\label{4.2}
\1^{\eta}_{D,E}=
\left\{
\begin{array}{ccl}
\1_{E}(k) & {\rm if} & \eta=(k)\\
\1_{D}(k) & {\rm if} & \eta=(k)*\\
\1_{D\cap E}(k)& {\rm if} & \eta=(k)\circ \\
1 & {\rm if} & \eta=(0)
\end{array}
\right.
\end{equation}
{\it for all $\eta \in {\cal T}$. Here, $\1_{A}$ denotes 
the indicator function of the set $A$.
Moreover, $I^{\eta}$ is a $D \cap E$-adapted regular process.}\\[5pt]
{\it Proof.} 
Denote $P=P^{(D)}$, $Q=P^{(E)}$. Since
$F\otimes G\in {\cal S}(D,E)$, 
and $A^{\eta}$ is ${\bf N}$-adapted, therefore, 
using the continuous tensor product
decomposition of exponential vectors, we obtain
$$
\langle
x,
F(s)\otimes G(s) \# (A_{t}^{\eta}-A_{s}^{\eta})
y
\rangle
$$
\begin{eqnarray*}
&=&
\langle 
x_{s]}, 
\widetilde{F}(s)\widetilde{G}(s)
y_{s]}
\rangle
\langle
x_{[s,t]}, 
P(A_{t}^{\eta}-A_{s}^{\eta}) Q
y_{[s,t]}
\rangle
\langle
x_{[t}, PQ
y_{[t}
\rangle \\
&=& 
\mu^{\eta}_{D,E}([s,t])
\langle 
x_{s]} \widetilde{F}(s)\widetilde{G}(s)
y_{s]}
\rangle
\langle
Px_{[s},
Qy_{[s}
\rangle\\
&=&\mu^{\eta}_{D,E}([s,t])
\langle
x, F(s)G(s)y
\rangle
\end{eqnarray*}
which completes the proof of the first part of the lemma.
The second part is obvious.
\hfill $\Box$\\
\indent{\par}
If we set $D=E=V={\bf N}$, then $\1_{D,E}^{\eta}\equiv 1$
for all $\eta \in {\cal T}$ and we obtain the formulas
of boson calculus for adapted processes (see [Mo-Si] or [P]).
Thus, the new feature in this lemma is that
apart from the time-dependent measures given by formula (\ref{4.1}),
which are the same as in boson calculus
and can be associated with time,
we also have the $0-1$ {\it color multipliers}
$\1^{\eta}_{D,E}$ given by formula (\ref{4.2}).
We choose to incorporate them into the measures $\mu_{D,E}^{\eta}$
(thus we can get trivial measures) for notational convenience. 

Let us now evaluate matrix elements of the form
$\langle I_{1}^{\eta_{1}}(t)x, I_{2}^{\eta_{2}}(t)y\rangle$
for $x=w\varepsilon (u)$, $y=z\varepsilon (v)$, $u,v\in {\cal H}_{0}$,
$w,z\in {\cal D}_{0}$, where
$$
I_{i}^{\eta_{i}}(t)=
\int_{0}^{t}F_{i}\otimes G_{i} \# dA^{\eta_{i}},
$$
for elementary simple biprocesses 
$F_{i}\otimes G_{i}\in {\cal S}(D_{i},E_{i})$, where
$D_{i},E_{i}\in {\cal P}({\bf N})$ and
$\eta_{i}\in {\cal T}$, $i=1,2$. 

Let us also use a short-hand notation
$$
\Delta_{1,2}=
\langle
P_{1}(A_{t}^{\eta_{1}}-A_{s}^{\eta_{1}})
Q_{1}x,
P_{2}(A_{t}^{\eta_{2}}-A_{s}^{\eta_{2}})
Q_{2}y
\rangle 
$$
where $P_{i}=P^{(D_{i})}$, $Q_{i}=P^{(E_{i})}$, $i=1,2$ and
$0\leq s < t <\infty$. 
We will also use the following convenient notation.
Instead of sets $E_{1}, D_{1}, D_{2}, E_{2}$
we will use integers $-2,-1,1,2$, respectively.
Moreover, we will identify $\1_{[-2]}\equiv \1_{E_{1}}$,
$\1_{[-1]}\equiv \1_{D_{1}}$, $\1_{[1]}\equiv \1_{D_{2}}$,
$\1_{[2]}\equiv \1_{E_{2}}$ and extend this notation
multiplicatively
$$
\1_{[n,m]}=\1_{[n]}\1_{[n+1]}\ldots \1_{[m]}
$$
for $n,m\in \{-2,-1,1,2\}$. 
Thus, for instance $\1_{D_{1}\cap D_{2}}=\1_{[-1,1]}$,
$\1_{D_{1}\cap D_{2}\cap E_{2}}=\1_{[-1,2]}$, etc.\\
\indent{\par}
{\sc Proposition 4.3.}
{\it Let $F_{i}\otimes G_{i}\in{\cal S}(D_{i},E_{i})$,
where $D_{i}, E_{i}\in {\cal P}({\bf N})$ 
and let $\eta_{i}\in {\cal T}$, where $i=1,2$. 
Then}
$$
\Delta_{1,2}=
[
\kappa_{1,2}
\mu^{\eta_{1}\dagger}_{E_{1},E_{2}}([s,t])
\mu^{\eta_{2}}_{E_{1},E_{2}}([s,t])
+
\mu_{1,2}([s,t])
] 
\langle 
P_{1}Q_{1}x,
P_{2}Q_{2}y
\rangle
$$
{\it where $\kappa_{1,2}\in \{0,1\}$ and 
the non-trivial part of the table of measures $\mu_{1,2}$ is given by}\\
\unitlength=1mm
\special{em:linewidth 0.4pt}
\linethickness{0.4pt}
\begin{picture}(110.00,30.00)(20.00,5.00)
\put(30.00,21.00){\line(1,0){120.00}}
\put(52.00,30.00){\line(0,-1){25.00}}
\put(35.00,25.00){$\mu_{1,2}$}
\put(65.00,25.00){$(k)*$}
\put(105.00,25.00){$(k)\circ$}
\put(35.00,15.00){$(k)*$}
\put(35.00,9.00){$ (k)\circ$}
\put(65.00,15.00)
{$\1_{[-1,1]}(k)\mu^{(0)}$}
\put(65.00,9.00)
{$\1_{[-2,1]}(k)
\mu^{(k)*}$}
\put(105.00,15.00)
{$\1_{[-1,2]}(k)
\mu^{(k)}$}
\put(105.00,9.00)
{$\1_{[-2,2]}(k)
\mu^{(k)\circ } $ }
\end{picture}\\[5pt]
{\it Proof.}
We only prove the case
$\eta_{1}=(k_{1})*$, $\eta_{2}=(k_{2})*$. Using
the relation
$$
P^{(V)}A_{t}^{(k)*}=\1_{V}(k)A_{t}^{(k)*}P^{(V)}
$$
which holds for all 
$V\in {\cal P}({\bf N})$, $k\in {\bf N}$ and $t\geq 0$, 
we obtain
\begin{eqnarray*}
\Delta_{1,2}
&=&
\1_{D_{1}}(k_{1})\1_{D_{2}}(k_{2})
\langle
(A_{t}^{(k_{1})*}-A_{s}^{(k_{1})*})P_{1}Q_{1}x,
(A_{t}^{(k_{2})*}-A_{s}^{(k_{2})*})P_{2}Q_{2}y
\rangle\\
&=&
[
\1_{E_{1}\cap D_{2}}(k_{2})
\1_{E_{2}\cap D_{1}}(k_{1})
\mu^{(k_{1})}([s,t])\mu^{(k_{2})*}([s,t])\\
&+&
\delta_{k_{1},k_{2}} \1_{D_{1}\cap D_{2}}(k_{1})(t-s)]
\langle
P_{1}Q_{1}x,
P_{2}Q_{2}y
\rangle\\
&=&
[\1_{D_{1}}(k_{1})\1_{D_{2}}(k_{2})
\mu^{\eta_{1}\dagger}_{E_{1},E_{2}}([s,t])\mu^{\eta_{2}}_{E_{1},E_{2}}([s,t])
+
\delta_{k_{1},k_{2}}
\1_{[-1,1]}(k_{1})(t-s)]\\
&\times&
\langle
P_{1}Q_{1}x,
P_{2}Q_{2}y
\rangle .
\end{eqnarray*}
The other cases are proved in a similar way. We do not give
the explicit formulas for $\kappa_{1,2}$ since they are not relevant
for the calculus (see Lemma 4.4). 
\hfill $\Box$\\
\indent{\par}
{\sc Lemma 4.4.} 
{\it Let $x=w\varepsilon(u)$, $y=z\varepsilon (v)$, where
$w,z\in {\cal D}_{0}$ and $u,v\in {\cal H}_{0}$.
Under the assumptions of Proposition 4.3 we have}
\begin{eqnarray*}
\langle 
I_{1}^{\eta_{1}}(t)x,
I_{2}^{\eta_{2}}(t)y
\rangle 
&=&
\int_{0}^{t}
\langle
F_{1}(s)G_{1}(s)x,
I_{2}^{\eta_{2}}(s) y
\rangle
d{\mu}_{1}(s)\\
&+&
\int_{0}^{t}
\langle
I_{1}^{\eta_{1}}(s) 
x,
F_{2}(s)G_{2}(s) 
y
\rangle
d{\mu}_{2}(s)\\
&+&
\int_{0}^{t}
\langle
F_{1}(s)G_{1}(s)
x,
F_{2}(s)G_{2}(s) 
y 
\rangle
d{\mu}_{1,2}(s)
\end{eqnarray*}
{\it where the measures $\mu_{1,2}$ are given by Proposition 4.3, and}
\begin{equation}
\label{4.3}
\mu_{1}=
\left\{
\begin{array}{ccl}
\1_{[-2]}(k_{1})\mu^{(k_{1})*} & {\rm if } & \eta_{1} =(k_{1})\\
\1_{[-1,2]}(k_{1})\mu^{(k_{1})} & {\rm if}& \eta_{1}=(k_{1})*\\
\1_{[-2,2]}(k_{1})\mu^{(k_{1})\circ} &{\rm if} & \eta_{1} =(k_{1})\circ\\
\mu^{(0)} &{\rm if} & \eta_{1} =(0)
\end{array}
\right.
\end{equation}
\begin{equation}
\label{4.4}
\mu_{2}=
\left\{
\begin{array}{ccl}
\1_{[2]}(k_{2})\mu^{(k_{2})} & {\rm if } & \eta_{2} =(k_{2})\\
\1_{[-2,1]}(k_{2})\mu^{(k_{2})*} & {\rm if}& \eta_{2}=(k_{2})*\\
\1_{[-2,2]}(k_{2})\mu^{(k_{2})\circ} &{\rm if} & \eta_{2} =(k_{2})\circ\\
\mu^{(0)} &{\rm if} & \eta_{2} =(0)
\end{array}
\right.
\end{equation}
{\it Proof.}
The proof is based on the Hudson-Parthasarathy theory
(for instance, see [H-P1] or [P]).
The main changes that come into play in our case are
due to 0-1 color multipliers which
produce more zeros in our formulas. 
We provide only the basic algebraic calculations.

Using standard decompositions of exponential vectors, we obtain
$$
x =
x_{t_{j-1}]}\otimes 
x_{[t_{j-1},t_{j}]}\otimes 
x_{[t_{j}}
$$
for each $j=1, \ldots ,n$, and an analogous formula for
$y$, which gives the formula
$$
\langle 
I_{1}^{\eta_{1}}(t)x, 
I_{2}^{\eta_{2}}(t)y
\rangle = S_{1} +S_{2} +S_{3}
$$
where
\begin{eqnarray*}
S_{1}
&=&
\sum_{j=1}^{n}
\langle 
F_{1}(t_{j-1})G_{1}(t_{j-1})x_{t_{j-1}]},
I_{2}^{\eta_{2}}(t_{j-1})y_{t_{j-1}]}
\rangle\\
&\times&
\langle
P_{1}(A^{\eta_{1}}(t_{j})-A^{\eta_{1}}(t_{j-1}))
Q_{1} x_{[t_{j-1},t_{j}]},
P_{2}Q_{2}
y_{[t_{j-1},t_{j}]}
\rangle\\
&\times&
\langle
P_{1}Q_{1}x_{[t_{j})}
P_{2}Q_{2}y_{[t_{j})}
\rangle .\\
S_{2}
&=&
\sum_{j=1}^{n}
\langle 
I^{\eta_{1}}(t_{j-1})
x_{t_{j-1}]},
F_{2}(t_{j-1})G_{2}(t_{j-1})
y_{t_{j-1}]}
\rangle\\
&\times&
\langle
P_{1}Q_{1}
x_{[t_{j-1},t_{j}]},
P_{2}
(A^{\eta_{2}}(t_{j})-A^{\eta_{2}}(t_{j-1}))
Q_{2} 
y_{[t_{j-1},t_{j}]}
\rangle\\
&\times&
\langle
P_{1}Q_{1}
x_{[t_{j})}
P_{2}Q_{2}
y_{[t_{j})}
\rangle .\\
S_{3}
&=&
\sum_{j=1}^{n}
\langle 
F_{1}(t_{j-1})G_{1}(t_{j-1})
x_{t_{j-1}]},
F_{2}(t_{j-1})G_{2}(t_{j-1})
y_{t_{j-1}]}
\rangle\\
&\times&
\langle
P_{1}(A^{\eta_{1}}(t_{j})-A^{\eta_{1}}(t_{j-1}))
Q_{1} x_{[t_{j-1},t_{j}]},\\
&&
P_{2}(A^{\eta_{2}}(t_{j})-A^{\eta_{2}}(t_{j-1}))
Q_{2} y_{[t_{j-1},t_{j}]}
\rangle\\
&\times&
\langle
P_{1}Q_{1}x_{[t_{j})}
P_{2}Q_{2}y_{[t_{j})}
\rangle .
\end{eqnarray*}
\indent{\par}
In each of those sums, the middle factor produces a 
complex-valued measure, denoted by $\mu_{1}, \mu_{2}$ and $\mu_{1,2}$,
respectively.
Note that Proposition 4.3  gives $\mu_{1,2}$. In turn,
$\mu_{1}$ and $\mu_{2}$ are defined by
\begin{eqnarray*}
\mu_{1}([s,t])
&=&
\langle
P_{1}(A^{\eta_{1}}(t)-A^{\eta_{1}}(s))
Q_{1} x_{[s,t]},
P_{2}Q_{2}
y_{[s,t]} 
\rangle\\
\mu_{2}([s,t])
&=&
\langle
P_{1}Q_{1}
x_{[s,t]},
P_{2}
(A^{\eta_{2}}(t)-A^{\eta_{2}}(s))
Q_{2} 
y_{[s,t]}
\rangle .
\end{eqnarray*}
It can be seen that it is enough to apply Lemma 4.2 to
get formulas (\ref{4.3})-(\ref{4.4}).
The remaining arguments are the same as in the usual case (see [P]).
\hfill $\Box$\\
\indent{\par}
{\it Remark.}
It is worth pointing out that the 0-1 color multipliers
which appear in all measures in both fundamental lemmas
follow easy-to-remember rules:\\[5pt]
\indent
{\it Rule 1.} all filters to the right of $dA^{(k)}$
must contain $k$,\\
\indent
{\it Rule 2.} all filters to the left of 
$dA^{(k)*}$ must contain $k$,\\
\indent
{\it Rule 3.} all filters on both sides of $dA^{(k)\circ}$ must contain $k$.\\[5pt]
We can of course assume that the integrated biprocesses satisfy these
rules and then we can skip the 0-1 color multipliers. However, when we
take a linear combination of biprocesses with different types of adaptedness
as integrands, it is important to keep track of what survives after
integration and what is killed, which is the main source
of non-commutativity leading, for instance, to $m$-free calculi.
\\[10pt]
\myownsection
\begin{center}
{\sc 5. An extension of the stochastic integral}
\end{center}
\indent{\par}
For given $D,E\in {\cal P}({\bf N})$ and $\eta \in {\cal T}$,
we will now take elementary 
simple biprocesses $X\in {\cal S}(D,E)$ as integrands 
and $A^{\eta}$ as integrators to approximate locally square integrable
$(D,E)$-adapted biprocesses in a topology specified below.
If $X=\sum_{i}F_{i}\otimes G_{i}\in {\cal A}(D,E)$ and 
the associated $D\cap E$-adapted process
$B=(B(t))_{t\geq 0}$ is given by 
$B(t)=\sum_{i}F_{i}(t)G_{i}(t)$
for all $t\geq 0$, then we will write $B\models X$.\\
\indent{\par}
{\sc Lemma 5.1.}
{\it Let $X\in {\cal S}(D,E)$, $x=w\varepsilon (u)$, 
where
$D,E\in {\cal P}({\bf N})$, $w\in {\cal D}_{0}$,
$u\in {\cal H}_{0}$, and let}
$$
I^{\eta }(t)x=
\int_{0}^{t}X\# dA^{\eta}x
$$ 
{\it where $\eta \in {\cal T}$. Then}
$$
\|I^{\eta}(t)x\|^{2}
\leq 
C^{\eta}_{D,E}(t)
\int_{0}^{t}
\|B(s)x\|^{2}\xi_{D,E}^{\eta}(ds)
$$
{\it where $B \models X$ and}
$$
\xi_{D,E}^{\eta}=|\sigma_{D,E}^{\eta}| +\nu_{D,E}^{\eta}
$$
{\it with $C^{\eta}_{D,E}(t)=e^{|\sigma_{D,E}^{\eta}|([0,t])}$,}
$$
\nu_{D,E}^{\eta}=\left\{
\begin{array}{cl}
\1_{D}(k)\mu^{(0)}& {\rm if}\;\eta =(k)*\\
\1_{D\cap E}(k)\mu^{(k)\circ} & {\rm if}\; \eta =(k)\circ \\
0 & {\rm otherwise}
\end{array}
\right.
$$
$$
\sigma_{D,E}^{\eta}=
\left\{
\begin{array}{ccl}
\1_{E}(k)\mu^{(k)} & {\rm if } & \eta =(k)\\
\1_{D\cap E}(k)\mu^{(k)*} & {\rm if}& \eta=(k)*\\
\1_{D\cap E}(k)\mu^{(k)\circ} &{\rm if} & \eta =(k)\circ\\
\mu^{(0)} &{\rm if} & \eta =(0)
\end{array}
\right.
$$
{\it and $|\mu|$ denoting the variation of $\mu$.}\\[5pt]
{\it Proof.}
This proof is based on Lemma 4.4 and is similar to that in
the adapted case (see [P]).\hfill $\Box$\\
\indent{\par}
{\sc Definition 5.2.}
For fixed $\eta\in {\cal T}$,
define a family of seminorms $\|.\|_{x,t,\eta}$ on ${\cal A}(D,E)$,
where $t\geq 0$, $x=w\varepsilon (u)\in {\cal M}_{0}$,
by
\begin{equation}
\label{5.1}
\|X\|_{x,t,\eta}^{2}=
\int_{0}^{t}\|B(s)x\|^{2}
\xi^{\eta}(ds)
\end{equation}
where $X\in {\cal A}(D,E)$ 
$t\geq 0$, $x\in {\cal M}_{0}$ and $\xi^{\eta}=\xi_{{\bf N},{\bf N}}^{\eta}$.
Denote by $L^{2}_{{\rm loc}}(D,E,dA^{\eta})$ the
linear vector space of all $(D,E)$-- adapted biprocesses
$X$ such that $\|X\|_{x,t,\eta}<\infty$
for all $x\in {\cal M}_{0}$ and $t\geq 0$.
We will say that $X$ is 
{\it locally square integrable with respect to} $dA^{\eta}$
(in our notation the dependence of this notion on the domain
$\widetilde{\cal E}_{0}$ and ${\bf R}^{+}$ is supressed).  \\
\indent{\par}
{\sc Theorem 5.3.}
{\it The stochastic integral with respect to the
fundamental process $A^{\eta}$ 
can be extended by continuity from 
$S(D,E)$ to $L^{2}_{{\rm loc}}(D,E,dA^{\eta})$
for any $\eta \in {\cal T}$ and $D,E \in {\cal P}({\bf N})$.}\\ [5pt]
{\it Proof.}
To show that ${\cal S}(D,E)$ is dense in
$L^{2}_{{\rm loc}}(D,E, dA^{\eta})$ 
for any $D,E \in {\cal P}({\bf N})$ and $\eta \in {\cal T}$, it is enough to
slightly modify the ideas of [Ac-Fa-Qu] and [H-P1], 
where we refer the reader for
details.\hfill $\Box$\\
\indent{\par}
For the approximating sequence of elementary simple biprocesses
$H^{(n)}\otimes K^{(n)}$ we set
$$
I^{\eta}(t)x= s-\lim_{n\rightarrow \infty}
I^{\eta}_{n}(t)x=
s-\lim_{n\rightarrow \infty}
\int_{0}^{t}H^{(n)}\otimes K^{(n)} \# dA^{\eta}
$$
on $\widetilde{\cal E}_{0}$. In view of Lemma 5.1, the sequence 
$(I_{n}^{\eta}(t)x)_{n\in {\bf N}}$ is Cauchy
for each $t\geq 0$ and $x\in {\cal M}_{0}$, hence convergent.
Moreover, the convergence is uniform for $t$ in finite intervals
and the limit does not depend on the choice of approximating simple
biprocesses. 
In this way we obtain
a $D\cap E$-adapted process $I^{\eta}$.\\
\indent{\par}
{\sc Theorem 5.4.}
{\it Lemmas 4.2, 4.4 and 5.1 remain true for locally square
integrable integrands.}\\[5pt]
{\it Proof.}
By Lemma 5.1, the convergence of $I_{n}^{\eta}$ to $I$ 
on $[0,T]$ is uniform for fixed $T$. Therefore
we conclude that Lemmas 4.2 and 4.4 as well as
the norm estimate of Lemma 5.1 hold for all 
biprocesses which are locally square integrable w.r.t. appropriate
fundamental processes in the sense of Definition 5.2. \hfill $\Box$
\\[10pt]
\myownsection
\begin{center}
{\sc 6. The filtered It$\hat{{\sc o}}$ formula}
\end{center}
In the multiplicative version of the filtered It$\hat{{\rm o}}$ formula 
we will use the differentials 
\begin{eqnarray*}
dI^{\eta_{1}}_{1}
&=&
G_{1}\otimes F_{1} \# dA^{\eta_{1}}=
G_{1}dA^{\eta_{1}}F_{1}\\
dI^{\eta_{2}}_{2}
&=&
F_{2}\otimes G_{2} \# dA^{\eta_{2}}=
F_{2}dA^{\eta_{2}}G_{2}
\end{eqnarray*}
instead of integrals $I^{\eta_{1}}_{1}$, $I^{\eta_{2}}_{2}$, respectively,
where $\eta_{i}\in {\cal T}$, $i=1,2$
(we switch the order of $F_{1}$ and $G_{1}$ in the first differential in order
to make a direct connection with Lemma 4.4).

Recall that the {\it boson It$\hat{o}$ table}
in the multivariate case is of the form\\
\unitlength=1mm
\special{em:linewidth 0.4pt}
\linethickness{0.4pt}
\begin{picture}(130.00,30.00)(-10.00,5.00)
\put(15.00,21.00){\line(1,0){90.00}}
\put(40.00,30.00){\line(0,-1){25.00}}
\put(18.00,25.00){$dA^{\eta_{1}}dA^{\eta_{2}}$}
\put(50.00,25.00){$dA^{(k)*}$}
\put(80.00,25.00){$dA^{(k)\circ}$}
\put(22.00,15.00){$dA^{(k)}$}
\put(22.00,9.00){$dA^{(k)\circ}$}
\put(50.00,15.00)
{$
dA^{(0)}
$}
\put(50.00,9.00)
{$dA^{(k)*}$}
\put(80.00,15.00)
{$
dA^{(k)}
$}
\put(80.00,9.00)
{$
dA^{(k)\circ}
$}
\end{picture}\\[5pt]
where we adopted the convention that only the 
non-trivial part of the It$\hat{{\rm o}}$ table is given,
thus the usual $\delta_{k_{1},k_{2}}$ does not appear here.
By $\bl A^{\eta_{1}}, A^{\eta_{2}} \br$ we will denote
the process obtained in the multiplication.

A pair of biprocesses, $(X,X^{\dagger})$, where 
$X\in {\cal A}(D,E)$ and $X^{\dagger}\in {\cal A}(E,D)$, will
be called an {\it adjoint pair} if 
$$
\langle x, B(s)y \rangle=
\langle B^{\dagger}(s)x, y\rangle
$$
for all $x,y\in {\cal M}_{0}$, where
$B\models X$ and $B^{\dagger}\models X^{\dagger}$ (cf. [P]).
A pair of $V$-adapted
processes, $(F,F^{\dagger})$ is an adjoint pair if 
$(F\otimes 1 , 1 \otimes F^{\dagger})$ is an adjoint
pair of biprocesses.
For instance, 
$F^{\dagger}=A^{(k,V)}, A^{(k,V)*}, A^{(k,V)\circ}$ or $A^{(0,V)}$ according to
whether $F=A^{(k,V)*}, A^{(k,V)},A^{(k,V)\circ}$ or $A^{(0,V)}$. Then
$(F,F^{\dagger})$ is an adjoint pair of $V$-adapted processes in $\Gamma$.
It follows easily from Lemma 4.2 that
if $X\in L^{2}_{{\rm loc}}(D,E, dA^{\eta})$,
$X^{\dagger}\in L^{2}_{{\rm loc}}(E,D,dA^{\eta \dagger})$ and
$(X,X^{\dagger})$ is an adjoint pair, then
$$
I_{1}(t)=\int_{0}^{t}X\# dA^{\eta},\;\;\;
I_{2}(t)=\int_{0}^{t}X^{\dagger}\# dA^{\eta \dagger}.
$$
is an adjoint pair of $D\cap E$-adapted processes.\\
\indent{\par}
{\sc Theorem 6.1.} ({\sc Filtered It$\hat{{\rm o}}$ Formula})
{\it Let 
$G_{1}\otimes F_{1}\in 
L^{2}_{{\rm loc}}(E_{1},D_{1},dA^{\eta_{1}})$,
$F_{1}^{*}\otimes G_{1}^{*} \in L^{2}_{{\rm loc}}(D_{1},E_{1},
dA^{\eta_{1} \dagger})$,
$F_{2}\otimes G_{2}\in 
L^{2}_{{\rm loc}}(D_{2},E_{2}, dA^{\eta_{2}})$,
where $D_{1},D_{2},E_{1},E_{2}\in {\cal P}({\bf N})$ and
$\eta_{1}, \eta_{2}\in {\cal T}$.
Suppose that $I^{\eta_{1}}_{1}I^{\eta_{2}}_{2}$ is a 
$D_{1}\cap D_{2}\cap E_{1}\cap E_{2}$- adapted process and that
$I^{\eta_{1}}_{1}F_{2}\otimes G_{2}$, 
$G_{1}\otimes F_{1}I^{\eta_{2}}_{2}$,
$G_{1}F_{1}F_{2} \otimes G_{2}$ and
$G_{1}\otimes F_{1}F_{2}G_{2}$ 
are locally square integrable with respect to $dA^{\eta_{2}}$, 
$dA^{\eta_{1}}$, $d\bl A^{\eta_{1}}, A^{\eta_{2}} \br$,
and $d\bl A^{\eta_{1}}, A^{\eta_{2}} \br$, respectively.
Then}
$$
d(I_{1}^{\eta_{1}} I_{2}^{\eta_{2}})
=
I_{1}^{\eta_{1}}dI_{2}^{\eta_{2}} +
dI_{1}^{\eta_{1}}I_{2}^{\eta_{2}} +
dI_{1}^{\eta_{1}}dI_{2}^{\eta_{2}}
$$
{\it where}
\begin{eqnarray*}
I_{1}^{\eta_{1}}dI_{2}^{\eta_{2}} 
&=&
I_{1}^{\eta_{1}}F_{2}\otimes G_{2} \# dA^{\eta_{2}}\\
dI_{1}^{\eta_{1}}I_{2}^{\eta_{2}} 
&=&
G_{1}\otimes F_{1}I_{2}^{\eta_{2}}\# dA^{\eta_{1}}
\end{eqnarray*}
{\it and the It$\hat{\rm o}$ correction can be written in two equivalent ways:}
\begin{eqnarray*}
dI_{1}^{\eta_{1}}dI_{2}^{\eta_{2}}
&=&
G_{1}\otimes \rho_{\eta_{1},\eta_{2}}(F_{1}F_{2})G_{2}\# 
d\bl A^{\eta_{1}}, A^{\eta_{2}}\br\\
&=&
G_{1}\rho_{\eta_{1},\eta_{2}}(F_{1}F_{2})\otimes G_{2}\# 
d\bl A^{\eta_{1}}, A^{\eta_{2}}\br
\end{eqnarray*}
{\it where 
$\rho_{\eta_{1}, \eta_{2}}(F_{1}F_{2})=\1_{[-1,1]}(k_{1})F_{1}F_{2}
\equiv\1_{D_{1}\cap D_{2}}(k_{1})F_{1}F_{1}$ 
for those values of $\eta_{1}, \eta_{2}$ which belong
to the non-trivial part of the It$\hat{{\rm o}}$ table, 
and for the other ones, it is zero.}\\[5pt]
{\it Proof.}
Without loss of generality assume that $I_{i}^{\eta_{i}}(0)=0$, $i=1,2$.
For all 
$t\in {\bf R}^{+}$, and
$x=w\varepsilon (u)$, $y=z\varepsilon (v)$, where
$w,y\in {\cal D}_{0}$ and $u,v\in {\cal H}_{0}$, we have
$$
\langle
x, I_{1}^{\eta_{1}}(t)I_{2}^{\eta_{2}}(t)y
\rangle
=\langle
I_{1}^{\eta_{1}\dagger}(t)x,
I_{2}^{\eta_{2}}(t)y
\rangle
$$
where
$$
I_{1}^{\eta_{1}\dagger}= 
\int_{0}^{t}F_{1}^{*}\otimes G_{1}^{*}
\# dA^{\eta_{1}\dagger}
$$
since $(I_{1}^{\eta_{1}}, I_{1}^{\eta_{1}\dagger})$
is an adjoint pair 
and the product $I_{1}^{\eta_{1}}I_{2}^{\eta_{2}}$
has $\widetilde{\cal E}_{0}$ in its domain.
By Lemma 4.4, this gives
\begin{eqnarray*}
\langle
x, I_{1}^{\eta_{1}}(t)I_{2}^{\eta_{2}}(t)y
\rangle
&=&\int_{0}^{t}
\langle 
I_{1}^{\eta_{1}\dagger}(s)x,
F_{2}(s)G_{2}(s)y
\rangle
d\mu_{2}(s)\\
&+&
\int_{0}^{t}
\langle
F_{1}^{*}(s)G_{1}^{*}(s)x,
I_{2}^{\eta_{2}}(s)y
\rangle 
d\mu_{1}(s)\\
&+&\int_{0}^{t}
\langle
F_{1}^{*}(s)G_{1}^{*}(s)
x,
F_{2}(s)G_{2}(s)
y
\rangle
d\mu_{1,2}(s)\\
&=&\int_{0}^{t}
\langle 
x,
I_{1}^{\eta_{1}}(s)F_{2}(s)G_{2}(s)y
\rangle
d\mu_{2}(s)\\
&+&
\int_{0}^{t}
\langle
x,
G_{1}(s)F_{1}(s)I_{2}^{\eta_{2}}(s)
y
\rangle 
d\mu_{1}(s)\\
&+&\int_{0}^{t}
\langle
x,
G_{1}(s)F_{1}(s)F_{2}(s)G_{2}(s)
y
\rangle
d\mu_{1,2}(s)
\end{eqnarray*}
where the measures $\mu_{1}$, $\mu_{2}$ and $\mu_{1,2}$ are 
determined appropriately depending on the fundamental 
processes $A^{\eta_{1}\dagger}$ and $A^{\eta_{2}}$ 
and are given by Lemma 4.4.
Let us look closer at the first two integrals.
Note that 
$$
F_{2}':=I_{1}^{\eta_{1}}F_{2}\in {\cal A}(D_{2}'),\;\;
F_{1}':=F_{1}I_{2}^{\eta_{2}}\in {\cal A}(D_{1}'),
$$
where $D_{2}'=D_{1}\cap D_{2}\cap E_{1}$ and 
$D_{1}'=D_{1}\cap D_{2}\cap E_{2}$.
Moreover, 
$$
F_{2}'\otimes G_{2}\in L^{2}_{{\rm loc}}(D_{2}', E_{2},
dA^{\eta_{2}}),\;\;{\rm and}\;\;
G_{1}\otimes F_{1}' \in L^{2}_{{\rm loc}}(E_{1}, D_{1}',
dA^{\eta_{1} \dagger})
$$ 
by assumption. Therefore,
we just need to check if the measures $\mu_{1}$, $\mu_{2}$ 
are also obtained when these biprocesses are integrated w.r.t. 
$A^{\eta_{2}}$ and $A^{\eta_{1}\dagger}$, respectively.
This is verified by using Lemma 4.2. For instance,
$$
\langle x,
\int_{0}^{t}F_{2}'\otimes G_{2}\# dA^{\eta_{2}}y
\rangle
=\int_{0}^{t}\langle 
x, F_{2}'(s)G_{2}(s)y
\rangle 
d\mu_{D_{2}', E_{2}}^{\eta_{2}}(s)
$$
where 
$$
\mu_{D_{2}',E_{2}}^{\eta_{2}}=
\mu_{D_{1}\cap D_{2} \cap E_{1}, E_{2}}^{\eta_{2}}=
\left\{
\begin{array}{lll}
\1_{[2]}(k_{2})\mu^{(k_{2})} & {\rm if} & \eta_{2}=(k_{2})\\
\1_{[-2,1]}(k_{2})\mu^{(k_{2})*} & {\rm if} &\eta_{2}=(k_{2})*\\
\1_{[-2,2]}(k_{2})\mu^{(k_{2})\circ} & {\rm if} 
&\eta_{2}=(k_{2})\circ\\
\mu^{(0)} & {\rm if} & \eta_{2}=(0)
\end{array}
\right.
$$ 
Comparing this with $\mu_{2}$ of Lemma 4.4, we 
conclude that $\mu_{2}=\mu_{D_{2}', E_{2}}^{\eta_{2}}$,
which enables us to write the first integral as
$$
\langle
x,
\int_{0}^{t} I_{1}^{\eta_{1}}F_{2}\otimes G_{2} \# dA^{\eta_{2}} 
y 
\rangle
$$
which, in differential notation, corresponds to 
$$
I_{1}^{\eta_{1}}dI_{2}^{\eta_{2}}=I_{1}^{\eta_{1}}F_{2}\otimes G_{2}\#
dA^{\eta_{2}}=I_{1}^{\eta_{1}}F_{2}dA^{\eta_{2}}G_{2}.
$$
A similar reasoning gives $\mu_{1}=\mu^{\eta_{1}}_{E_{1},D_{1}'}$,
which enables us to write the differential of the second integral 
as  
$$
dI_{1}^{\eta_{1}}I_{2}^{\eta_{2}}=G_{1}\otimes F_{1}I_{2}^{\eta_{2}}\#
dA^{\eta_{1}}=G_{1}dA^{\eta_{1}}F_{1}I_{2}^{\eta_{2}}.
$$
Let us finally evaluate the It$\hat{{\rm o}}$ correction.
This boils down to straightforward 
examination of 4 non-trivial 
cases. For instance,
let $\eta_{1}=(k_{1})$ and $\eta_{2}=(k_{2})\circ$.
Then 
$$
G_{1}\otimes F_{1}F_{2}G_{2}\in L^{2}_{{\rm loc}}
(E_{1}, D_{1}\cap D_{2}\cap E_{2}, dA^{(k_{1})})
$$
by assumption 
since $\bl A^{(k_{1})}, A^{(k_{2})\circ} \br=\delta_{k_{1},k_{2}}
A^{(k_{1})}$. We have
$$
\mu_{1,2}
=
\delta_{k_{1},k_{2}}\1_{[-1,2]}(k_{1})
\mu^{(k_{1})}
=
\delta_{k_{1},k_{2}}
\1_{[-1,1]}(k_{1})
\mu_{E_{1}, D_{1}\cap D_{2}\cap E_{2}}^{(k_{1})}
$$
This gives
$$
dI_{1}^{(k_{1})}dI_{2}^{(k_{2})\circ}
=
\delta_{k_{1},k_{2}}\1_{[-1,1]}(k_{1})
G_{1}\otimes F_{1}F_{2}G_{2} \#
dA^{(k_{1})}.
$$
The second formula for the Ito correction as well as other cases
are proved in a similar way.
\hfill $\Box$\\
\indent{\par}
{\sc Definition 6.2.}
A $(D,E)$-adapted stochastic biprocess $X=(X(t))_{t\geq 0}$ 
will be called {\it locally-bounded} if
\indent{\par}
(1) $X(t)\in ({\cal B}(h_{0})\otimes {\cal B}({\Gamma}))
\otimes ({\cal B}(h_{0})\otimes {\cal B}({\Gamma}))$
for all $t\geq 0$
\indent{\par}
(2) the map $t\rightarrow B(t)x$ is strongly measurable
for all $x\in \widetilde{\Gamma}({\cal H})$, where $B\models X$
\indent{\par}
(3) $\sup_{0\leq s \leq t}\|B(s)\|<\infty$ for all $t>0$.
\indent{\par}
The space of locally bounded $(D,E)$-adapted
biprocesses will be denoted by ${\cal B}_{{\rm loc}}(D,E)$.\\
\indent{\par}
{\it Remark.}
One can give a shorter formulation of Theorem 6.1 if one makes stronger
assumptions. Namely, if we assume that
$G_{1}\otimes F_{1}\in {\cal B}_{{\rm loc}}(E_{1},D_{1})$
and $F_{2}\otimes G_{2}\in {\cal B}_{{\rm loc}}(D_{2},E_{2})$, 
then the assumptions of Theorem
6.1 are satisfied.

Let us now give conditions under which an infinite
sum of stochastic integrals associated with the same pair of filters
$(D,E)$ but different $\eta \in {\cal T}$ is well-defined on
the domain $\widetilde{\cal E}_{0}$. Namely, we want to define
$$
I(t)=
\sum_{\eta\in {\cal T}}
\int_{0}^{t}X^{\eta}\# dA^{\eta}
$$
on $\widetilde{\cal E}_{0}$
for all $t\in {\bf R}^{+}$, where
$X^{\eta}$ is $(D,E)$- adapted and locally square integrable with
respect to $dA^{\eta}$ for each $\eta\in {\cal T}$ and the sum
is possibly infinite.
The case $D=E={\bf N}$ was studied in [Mo-Si] (see also [P]).

An approximating sequence of integrals will be given by
$$
I^{n}(t)=\sum_{\eta \in {\cal T}(n)}\int_{0}^{t}X^{\eta}\#dA^{\eta}
$$
where 
$$
{\cal T}(n)=\{(k),(l)*,(r)\circ , (0): 1\leq k,l,r \leq n\}
$$
denotes the set of indices associated with fundamental processes
of colors less than or equal to $n$. Now, let
$$
N(u)={\rm max}
\{
k: u^{(k)}\;\; {\rm is\; a \; non-zero\; function
\; in}\; L^{2}({\bf R}^{+})\}
$$
with
$$
{\cal T}(n,u)=\{(k), (l)*, (r)\circ , (0): 1\leq k,r \leq n\wedge N(u),
1\leq l \leq n\}
$$
for $u\in {\cal H}_{0}$, and we set ${\cal T}(u)={\cal T}(\infty, u)$.
Note that in the case of creation process there is no constraint
on the color support of $u$ -- the reason is that the 
It$\hat{{\rm o}}$ correction term corresponding
to the creation-- creation pair always appears in 
the filtered It$\hat{{\rm o}}$ formula
irrespective of $u$.\\
\indent{\par}
{\sc Theorem 6.3.}
{\it Suppose $X^{\eta}\in L^{2}_{{\rm loc}}(D,E,dA^{\eta})$
for each $\eta\in {\cal T}$ and that}
\begin{equation}
\label{6.1}
\sum_{\eta\in {\cal T}}
\int_{0}^{t}
\|B^{\eta}(s)x\|^{2}d\nu_{u}(s) <\infty.
\end{equation}
{\it for all $t\geq 0 $ and $x\in {\cal M}_{0}$,
where the real-valued measures $\nu_{u}$, $u\in {\cal H}_{0}$, are given by}
$$
\nu_{u}(t)=\int_{0}^{t}(\sum_{k\in {\bf N}}|u^{(k)}(s)|^{2}+1)ds
$$
{\it and $B^{\eta} \models X^{\eta}$.
Then there exists a regular $D\cap E$-- adapted process
$I$ such that}
$$
\lim_{n\rightarrow \infty}{\rm sup}_{0\leq s \leq t}
\|I^{(n)}(s)x -I(t)x\|=0
$$
$$
\|I(t)x \|^{2}
\leq 2e^{\nu_{u}(t)}
\sum_{\eta \in {\cal T}(u)} 
\int_{0}^{t}
\|B^{\eta}(s)x\|^{2}d\nu_{u}(s)
$$
{\it for all $x\in {\cal M}_{0}$ and $t\geq 0$.}\\[5pt]
{\it Proof.}
The proof is similar to that in [P]. \hfill $\Box$\\[10pt]
\myownsection
\begin{center}
{\sc 7. The $m$-free calculi}
\end{center}
In order to include $m$-free calculi for $1\leq m \leq \infty$, 
we need to add 
biprocesses associated with different pairs of filters $(D,E)$ and 
integrate them w.r.t. fundamental processes. In this section we use
filtered stochastic calculus to recover $m$-free calculi.\\
\indent{\par}
{\sc Definition 7.1.}
For all $m\in {\bf N}^{*}$ and $\alpha \in {\cal F}_{m}$
we define on $\widetilde{\cal E}_{0}$
the integrals w.r.t. extended
$m$-free fundamental processes by the linear extension of
\begin{eqnarray*}
\int_{0}^{t}X\# dl^{\alpha}:&=&
\sum_{(\eta , V)\sim \alpha}
\int_{0}^{t}X[\eta, V]\# dA^{\eta}
\end{eqnarray*}
where $X=F\otimes G$ is a $(D,E)$-adapted stochastic biprocess
for which the integrands on the RHS are 
locally square integrable biprocesses and $X[\eta ,V]$
is given by (\ref{3.2}).
We will say that $X$ is
{\it locally square integrable w.r.t.} 
$dl^{\alpha}$.
The space spanned by such biprocesses will be denoted by
$L^{2}_{{\rm loc}}(D,E,dl^{\alpha})$. Note that if
$X$ is locally bounded, then 
$X\in L^{2}_{{\rm loc}}(D,E,dl^{\alpha})$ for all $\alpha\in {\cal F}_{m}$.
It can be shown that the integrals w.r.t. extended {\it free} 
fundamental processes always reduce to finite sums on $\widetilde{\cal E}_{0}$
and thus are well-defined.\\
\indent{\par}
{\sc Proposition 7.2.}
{\it Let $x=w\varepsilon (u)$, $y=z\varepsilon (v)$, where
$w,z\in {\cal D}_{0}$, 
$u,v\in {\cal H}_{0}$, $\alpha \in {\cal F}_{m}$ and
$m\in {\bf N}^{*}$, and assume that 
$F\otimes G\in L^{2}_{{\rm loc}}(D,E, dl^{\alpha})$,
where $D,E \in {\cal P}({\bf N})$. Then}
$$
\langle 
x , \int_{0}^{t}F\otimes G \# dl^{\alpha} 
y
\rangle
=
\int_{0}^{t}
(
x,
F(s)\otimes G(s)
y
)
\# d\widehat{\nu}^{\alpha}(s)
$$
{\it where}
\begin{displaymath}
\widehat{\nu}^{\alpha}
=
\left\{
\begin{array}{lll}
\sum_{k\in E(m)}
\dl I, P^{[k-1]} \dr \mu^{(k)} &{\rm if} & \alpha =(m)\\
\sum_{k\in D(m)}
\dl I, P^{[k-1]} \dr \mu^{(k)*} & {\rm if} & \alpha =(m)*\\
\sum_{k \in D(m)\cap E(m)}
\dl I, P^{[k]} \dr \mu^{(k)\circ} & {\rm if} & \alpha =(m)\circ\\
\dl I ,P^{(m)}\dr \mu^{(0)} & {\rm if} & \alpha =(m)\cdot
\end{array}
\right.
\end{displaymath}
{\it and $D(m)=D\cap \{1, \ldots , m\}$}.\\[5pt]
{\it Proof.}
It is a straightforward consequence of (2.6)-(2.9) and Lemma 4.2.
\hfill $\Box$\\
\indent{\par}
Before we state the general version of the $m$-free It$\hat{{\rm o}}$ formula,
let us first establish its easy case, 
the {\it $m$-free It$\hat{{\rm o}}$ table}\\
\indent{\par}
{\sc Proposition 7.3.}
{\it Let $\alpha_{1}, \alpha_{2}\in {\cal F}_{m}$, where
$m\in {\bf N}^{*}$. Then
$l^{\alpha_{1}}l^{\alpha_{2}}$ satisfies the relation}
$$
d(l^{\alpha_{1}} l^{\alpha_{2}})=
l^{\alpha_{1}}dl^{\alpha_{2}} +
dl^{\alpha_{1}}l^{\alpha_{2}}
+
dl^{\alpha_{1}}dl^{\alpha_{2}} 
$$
{\it where the It$\hat{o}$ correction is given by the following
multiplication table:}\\
\unitlength=1mm
\special{em:linewidth 0.4pt}
\linethickness{0.4pt}
\begin{picture}(80.00,30.00)(-27.00,5.00)
\put(10.00,21.00){\line(1,0){75.00}}
\put(32.00,30.00){\line(0,-1){25.00}}
\put(13.00,25.00){$dl^{\alpha_{1}} dl^{\alpha_{2}}$}
\put(45.00,25.00){$dl^{(m)*}$}
\put(65.00,25.00){$dl^{(m)\circ}$}
\put(18.00,15.00){$dl^{(m)}$}
\put(18.00,9.00){$dl^{(m)\circ}$}
\put(45.00,15.00)
{$dl^{(m)\cdot}$}
\put(45.00,9.00)
{$dl^{(m)*}$}
\put(65.00,15.00)
{$dl^{(m)}$}
\put(65.00,9.00)
{$dl^{(m)\circ}$}
\end{picture}\\[5pt]
{\it with the convention that only 
the non-trivial part of the table is given.}\\[5pt]
{\it Proof.}
The proof is based on the It$\hat{{\rm o}}$ table for boson 
calculus (see Section 7). 
Assume first that $m$ is finite.
We will consider one case,
leaving the other ones to the reader. Using bi-linearity, we obtain
\begin{eqnarray*}
d(l^{(m)}l^{(m)*})
&=&
\sum_{k,l=1}^{m}P^{[k-1]}d(A^{(k)}A^{(l)*})P^{[l-1]}\\
&&
\sum_{k,l=1}^{m}
\left\{
P^{[k-1]}A^{(k)}dA^{(l)*}P^{[l-1]}
+
P^{[k-1]}dA^{(k)}A^{(l)*}P^{[l-1]}
\right.\\
&+&
\left.
P^{[k-1]}dA^{(k)}dA^{(l)*}P^{[l-1]}
\right\}\\
&=&
l^{(m)}dl^{(m)*}+dl^{(m)}l^{(m)*}
+ dA^{(0)}
\sum_{k=1}^{m}P^{[k-1]}\\
&=&
l^{(m)}dl^{(m)*}+dl^{(m)}l^{(m)*}
+ dl^{(m)\cdot}
\end{eqnarray*}
where we used the fact that the time differential $dA^{(0)}\equiv dt$ 
commutes with the projections $P^{[k-1]}$ for all $k\in {\bf N}$.
This gives $dl^{(m)}dl^{(m)*}=dl^{(m)\cdot}$. The strong
limit as $m\rightarrow \infty$ on $\widetilde{\cal E}_{0}$ 
of this relation gives $dl dl^{*}=dl^{\cdot}\equiv dt$.
The other cases are analogous.\hfill $\Box$\\
\indent{\par}
The processes obtained from the above It$\hat{{\rm o}}$ table will be denoted
$\bl l^{\alpha_{1}}, l^{\alpha_{2}} \br$, as in Section 6 for filtered
fundamental processes, i.e. 
$$
dl^{\alpha_{1}}dl^{\alpha_{2}}=d\bl l^{\alpha_{1}} , l^{\alpha_{2}} \br.
$$
Another notation will also be needed in the It$\hat{{\rm o}}$ formula.
Namely, for a given densely defined $V$-adapted 
linear operator $H$ on $\widetilde{\Gamma}({\cal H})$,
let
\begin{eqnarray}
\label{7.1}
\tr (H)
&=&
\sum\limits_{k\in V}P^{[k-1]}HP^{[k-1]}\\
\label{7.2}
\trr (H)
&=&
\sum\limits_{k\in V}P^{[k]}HP^{[k]}.
\end{eqnarray}
The two operators $\tr$ and
$\trr$ play a role of ``partial traces'' .
Of course, the only difference between $\tr$ and $\trr$ is that
the first one includes the vacuum in its trace if $k=1\in V$, whereas
$\trr$ doesn't. 

The differentials which enter the $m$-free It$\hat{{\rm o}}$ formula are given by
$$
dJ_{1}^{\alpha_{1}}=G_{1}\otimes F_{1}\# dl^{\alpha_{1}}=
G_{1}dl^{\alpha_{1}}F_{1}
$$
$$
dJ_{2}^{\alpha_{2}}=F_{2}\otimes G_{2}\# dl^{\alpha_{2}}
=F_{2}dl^{\alpha_{2}}G_{2}
$$
where $\alpha_{1},\alpha_{2}\in {\cal F}_{m}$. They correspond
to integrals $J_{1}^{\alpha_{1}}$, $J_{2}^{\alpha_{2}}$, respectively, for
which we assume, without loss of generality, that $J_{1}^{\alpha_{1}}(0)=
J_{2}^{\alpha_{2}}(0)=0$.
\\
\indent{\par}
{\sc Theorem 7.4.}
{\it Let 
$G_{1}\otimes F_{1}
\in L^{2}_{{\rm loc}}(E_{1},D_{1},dl^{\alpha_{1}})$,
$F_{1}^{*}\otimes G_{1}^{*}
\in L^{2}_{{\rm loc}}(D_{1},E_{1},dl^{\alpha_{1}\dagger})$,
$F_{2}\otimes G_{2}\in L^{2}_{{\rm loc}}(D_{2},E_{2}, dl^{\alpha_{2}})$,
where $D_{1}$,$D_{2}$,$E_{1}$,$E_{2}$$\in {\cal P}({\bf N})$ and 
$\alpha_{1}, \alpha_{2}\in {\cal F}_{m}$, $m\in {\bf N}^{*}$. 
Suppose that $J_{1}^{\alpha_{1}}J_{2}^{\alpha_{2}}$ is a filtered
adapted process and that $J^{\alpha_{1}}F_{2}\otimes G_{2}$,
$G_{1}\otimes F_{1}J_{2}^{\alpha_{2}}$, $G_{1}F_{1}F_{2}\otimes G_{2}$
and $G_{1}\otimes F_{1}F_{2}G_{2}$ are locally square integrable
with respect to $dl^{\alpha_{2}}$, $dl^{\alpha_{1}}$,
$d\bl l^{\alpha_{1}} , l^{\alpha_{2}} \br $ and 
$d\bl l^{\alpha_{1}} , l^{\alpha_{2}} \br $, respectively. Then}
$$
d(J_{1}^{\alpha_{1}}J_{2}^{\alpha_{2}})
=J_{1}^{\alpha_{1}}dJ_{2}^{\alpha_{2}}
+ dJ_{1}^{\alpha_{1}}J_{2}^{\alpha_{2}}
+ dJ_{1}^{\alpha_{1}}dJ_{2}^{\alpha_{2}}
$$
{\it where}
\begin{eqnarray*}
dJ_{1}^{\alpha_{1}}J_{2}^{\alpha_{2}}
&=&
G_{1}\otimes F_{1}J_{2}^{\alpha_{2}}\# dl^{\alpha_{1}},\\
J_{1}^{\alpha_{1}}dJ_{2}^{\alpha_{2}}
&=&
J_{1}^{\alpha_{1}}F_{2}\otimes G_{2} \# dl^{\alpha_{2}},
\end{eqnarray*}
{\it and the It$\hat{o}$ correction can be written
in two equivalent ways:}
\begin{eqnarray*}
dJ_{1}^{\alpha_{1}}dJ_{2}^{\alpha_{2}}
&=&
G_{1}\otimes P_{\alpha_{1},\alpha_{2}}(F_{1}F_{2})G_{2}
\# d\bl l^{\alpha_{1}}, l^{\alpha_{2}} \br\\
&=&
G_{1}P_{\alpha_{1},\alpha_{2}}(F_{1}F_{2})\otimes G_{2}
\# d\bl l^{\alpha_{1}}, l^{\alpha_{2}} \br 
\end{eqnarray*}
{\it where}\\
\unitlength=1mm
\special{em:linewidth 0.4pt}
\linethickness{0.4pt}
\begin{picture}(80.00,30.00)(-27.00,5.00)
\put(10.00,21.00){\line(1,0){75.00}}
\put(32.00,30.00){\line(0,-1){25.00}}
\put(13.00,25.00){$P_{\alpha_{1},\alpha_{2}}$}
\put(45.00,25.00){$(m)*$}
\put(65.00,25.00){$(m)\circ$}
\put(13.00,15.00){$(m)$}
\put(13.00,9.00){$(m)\circ$}
\put(45.00,15.00)
{$\tr$}
\put(45.00,9.00)
{$\trr$}
\put(65.00,15.00)
{$\trr$}
\put(65.00,9.00)
{$\trr$}
\end{picture}\\[5pt]
{\it Proof.}
We will restrict our attention to the It$\hat{{\rm o}}$ correction.
Let $m$ be finite
and take $\alpha_{1}=(m)$ and $\alpha_{2}=(m)*$. Using
bi-linearity and Theorem 6.1, we obtain
\begin{eqnarray*}
dJ_{1}^{(m)}dJ_{2}^{(m)*}
&=&
G_{1}\otimes F_{1}F_{2}\otimes G_{2}\#
\dl 
\sum_{k=1}^{m}P^{[k-1]}dA^{(k)},
\sum_{l=1}^{m}dA^{(l)*}P^{[l-1]}
\dr\\
&=&
\sum_{k,l=1}^{m}
G_{1}P^{[k-1]}\otimes F_{1}F_{2} \otimes P^{[l-1]}G_{2}
\#
\dl
dA^{(k)}, dA^{(l)*}
\dr\\
&=&
\sum_{k\in D_{1}(m)\cap D_{2}(m)}
G_{1}P^{[k-1]}\otimes F_{1}F_{2}P^{[k-1]}G_{2}\# dA^{(0)}\\
&=&
\sum_{k\in D_{1}\cap D_{2}} 
G_{1}\otimes P^{[k-1]}F_{1}F_{2}P^{[k-1]}G_{2}\# dl^{(m)\cdot}\\
&=&
G_{1}\otimes \tr (F_{1}F_{2}) G_{2}\# dl^{(m)\cdot}.
\end{eqnarray*}
Using the second way of writing the It$\hat{{\rm o}}$ correction as 
indicated in Theorem 6.1 we obtain a similar expression in which
the differential $dl^{(m)\cdot}$ is placed right before $G_{2}$.

Let now $\alpha_{1}=(m)\circ$ and $\alpha_{2}=(m)*$.
Then
\begin{eqnarray*}
dJ_{1}^{(m)\circ}dJ_{2}^{(m)*}
&=&
G_{1}\otimes F_{1}F_{2}\otimes G_{2}\#
\dl
\sum_{k=1}^{m}dA^{(k)\circ}P^{[k]},
\sum_{l=1}^{m}dA^{(l)*}P^{[l-1]} 
\dr\\
&=&
\sum_{k,l=1}^{m}
G_{1}\otimes P^{[k]}F_{1}F_{2}\otimes P^{[k-1]}G_{2}\#
\dl dA^{(k)\circ}, dA^{(l)*} \dr\\
&=&
\sum_{k\in D_{1}(m)\cap D_{2}(m)}
G_{1}P^{[k]}F_{1}F_{2}\otimes P^{[k-1]}G_{2}\# dA^{(k)*}\\
&=&
\sum_{k\in D_{1}(m)\cap D_{2}(m)}
G_{1}P^{[k]}F_{1}F_{2}P^{[k]}\otimes G_{2}\# dA^{(k)*}P^{[k-1]}
\end{eqnarray*}
where we used the fact that for any $x\in \widetilde{\Gamma}({\cal H})$,
we have $x'=P^{[k-1]}x \in \widetilde{\Gamma}({\cal H}^{(k)})$ 
and thus 
$$
x''=dA^{(k)*}x'\in \widetilde{\Gamma}({\cal H}^{(k+1)})\ominus
\widetilde{\Gamma}({\cal H}^{(k)})
$$ 
since
$dA^{(k)*}$ adds color $k$ to the given vector and thus $P^{[k]}x''=x''$. 
Therefore, we obtain
\begin{eqnarray*}
dJ_{1}^{(m)\circ}dJ_{2}^{(m)*}
&=&
\sum_{k\in D_{1}\cap D_{2}}
G_{1}P^{[k]}F_{1}F_{2}P^{[k]}\otimes G_{2} 
\# (\sum_{r=1}^{m}dA^{(r)*}P^{[r-1]})\\
&=&
G_{1}\trr (F_{1}F_{2})\otimes G_{2} \# dl^{(m)*}
\end{eqnarray*}
which finishes the proof if the differential $dl^{(m)*}$
is to be written right before $G_{2}$. The second
formula for this product as well as the two
remaining cases are proved in a similar manner.
If $m=\infty$,
one uses Definition 3.1.(assumption 3) and
the fact that $m$-free processes converge strongly to free processes.
\hfill $\Box$\\[10pt]
\myownsection
\begin{center}
{\sc 8. Stochastic differential equations}
\end{center}
Let ${\cal P}_{0}$ be a finite collection of subsets of
the power set ${\cal P}({\bf N})$ which is closed under intersections.
We will establish existence and uniqueness of solutions of 
systems of stochastic differential equations of the type
\begin{eqnarray}
\label{8.1}
dI_{V}
&=& 
\sum_{C\cap D\cap E= V}
\sum_{\eta \in {\cal T}}
X^{\eta}_{C,D}I_{E}\# dA^{\eta}\\
\label{8.2}
I_{V}(0)
&=&
I_{V}^{(0)}
\end{eqnarray}
on $\widetilde{\cal E}_{0}$,
where $V\in {\cal P}_{0}$,
$X^{\eta}_{C,D}$ are suitable $(C,D)$-adapted 
locally bounded biprocesses and 
\begin{equation}
\label{8.3}
I_{V}^{(0)}=
\bar{I}_{V}^{(0)}\otimes P^{(V)}\in {\cal B}(h_{0})\otimes {\cal B}({\Gamma})
\end{equation}
where $\bar{I}_{V}^{(0)}$ are bounded operators on $h_{0}$.
By $\sum_{C\cap D\cap E = V}$ we understand the summation
over all $C,D,E\in {\cal P}_{0}$ such that $C\cap D \cap E = V$.
Thus, we suppress in the notation the fact that 
$C,D,E\in {\cal P}_{0}$. This convention will also be adopted in 
expressions of similar type in the sequel.
Of course, the above system of 
stochastic differential equations should be interpreted as the
system of stochastic integral equations
\begin{equation}
\label{8.4}
I_{V}(t)
= I_{V}^{(0)} +
\sum_{C\cap D\cap E = V}
\sum_{\eta \in {\cal T}}
\int_{0}^{t}
X^{\eta}_{C,D}I_{E}\# dA^{\eta}
\end{equation}
on $\widetilde{\cal E}_{0}$, where $V\in {\cal P}_{0}$ and $t\geq 0$. 
Note that on both sides of (8.1) and (8.4) we have 
processes of the same type of adaptedness. The reason for this is
that in a variety of interesting cases, processes of 
different types of adaptedeness are linearly independent
(see Lemma 8.3).

We will need the elementary estimate given below.\\
\indent{\par}
{\sc Proposition 8.1.}
{\it Let $B(t), L(t)\in {\cal B}(h)$
for all $t\geq 0$, where $h$ is a separable Hilbert space 
and let $\mu$ be a numerically-valued measure on ${\bf R}^{+}$ of bounded
variation. Suppose that the mapping $t\rightarrow B(t)y$ 
is strongly measurable for all $y\in h$ and locally bounded and
that the mapping $t\rightarrow L(t)x$ is strongly measurable, where
$x\in h$.
Then}
$$
\int_{0}^{t}\|B(s)L(s)x\|^{2}d\mu(s)
\leq 
\sup_{0\leq r \leq t}\|B(r)\|^{2}\int_{0}^{t}
\|L(s)x\|^{2}d\mu(s)
$$
{\it for all $t\geq 0$ .}\\[5pt]
{\it Proof.} Obvious.
\hfill $\Box$\\
\indent{\par}
{\sc Theorem 8.2.}
{\it Let $X_{C,D}^{\eta}\in {\cal B}_{{\rm loc}}(C,D)$ 
for all $\eta \in {\cal T}$ and $C,D\in {\cal P}_{0}$ and
let $I^{(0)}_{V}$ be given by formula (\ref{8.3}) for all
$V\in {\cal P}_{0}$. If} 
\begin{equation}
\label{8.5}
\sum_{\eta \in {\cal T}}\sup_{0\leq s \leq t}
\|B_{C,D}^{\eta}(s)\|^{2}<\infty
\end{equation}
{\it for all $t\geq 0 $ and $C,D\in {\cal P}_{0}$, where
$B_{C,D}^{\eta}\models X_{C,D}^{\eta}$, then
there exists a unique family of $V$-adapted
regular processes 
$(I_{V})_{V\in {\cal P}_{0}}$ satisfying equation (\ref{8.4}).}\\[5pt]
{\it Proof.}
The iterative scheme is established in the usual way. Thus, let
$I_{V}^{(0)}$ be the zero-th order approximation of 
$I_{V}$, $V\in {\cal P}_{0}$, and let
\begin{equation}
\label{8.6}
I_{V}^{(n)}(t)=I_{V}^{(0)}+
\sum_{C\cap D\cap E= V}\sum_{\eta \in {\cal T}}
\int_{0}^{t}
X_{C,D}^{\eta}I_{E}^{(n-1)}\# dA^{\eta}
\end{equation}
for all $n\in {\bf N}$. 
We will show that the iterative scheme is well-defined,
each $I_{V}^{(n)}$ is a regular $V$-adapted process 
and that the following estimate holds:
\begin{equation}
\label{8.7}
\|(I_{V}^{(n)}(t)-I_{V}^{(n-1)}(t))x\|^{2}
\leq 2^{n}|{\cal P}_{0}|^{6n}k_{T}^{n}l_{0}
\|x\|^{2}
\frac{(\nu_{u}(t))^{n}}{n!}
\end{equation}
for each $n\in {\bf N}$, $t\in [0,T]$ and $x\in {\cal M}_{0}$, 
where $|{\cal P}_{0}|$ denotes the cardinality of ${\cal P}_{0}$,
\begin{eqnarray*}
l_{0}
&=&
\max_{E\in {\cal P}_{0}} \|I_{E}^{(0)}\|^{2}\\
k_{T}
&=&
\max_{C,D\in {\cal P}_{0}}k_{T}(C,D)\\
k_{T}(C,D)
&=&
\sum_{\eta \in {\cal T}}\sup_{0\leq t \leq T}
\|B_{C,D}^{\eta}(t)\|^{2} \\
\end{eqnarray*}
and $\nu_{u}(t)=\nu_{u}([0,t])$. 

First of all, $X_{C,D}^{\eta}I_{E}^{(0)}$ is locally square integrable
w.r.t. $dA^{\eta}$, $\eta \in {\cal T}$, for all
$C,D,E\in {\cal P}_{0}$ since $X_{C,D}^{\eta}\in {\cal B}_{{\rm loc}}(D,E)$
and $I_{E}^{(0)}$ is a bounded operator
(Proposition 8.1 is used to show that the seminorms given by 
(\ref{5.1}) are finite).
Therefore,
$I_{V}^{(1)}$ will be well-defined if the estimate (\ref{6.1})
is proven to hold. This will follow from the proof
of the estimate (\ref{8.7}) for $n=1$ given below.
For fixed $T\geq 0$ and $0\leq t \leq T$ we obtain
$$
\|(I_{V}^{(1)}(t)-I_{V}^{(0)}(t))x\|^{2}
$$
\begin{eqnarray*}
&\leq &
|{\cal P}_{0}|^{6}\max_{C\cap D\cap E= V}
\|
\sum_{\eta \in {\cal T}}
\int_{0}^{t}X_{C,D}^{\eta}I_{E}^{(0)}\#
dA^{\eta}x
\|^{2}\\
&\leq &
2|{\cal P}_{0}|^{6}e^{\nu_{u}(T)}
\max_{C\cap D\cap E= V}
\left\{
\sum_{\eta \in {\cal T}}
\sup_{0\leq s \leq t} 
\|B_{C,D}^{\eta}(s)\|^{2}\|I_{E}^{(0)}\|^{2}
\right\}
\|x\|^{2}\nu_{u}(t)\\
&\leq &
2|{\cal P}_{0}|^{6}e^{\nu_{u}(T)}
k_{T}l_{0}
\|x\|^{2}
\nu_{u}(t)
\end{eqnarray*}
using the estimate (\ref{6.1}) and Proposition 8.1,
which proves the estimate for $n=1$. Moreover, note that 
$I_{V}^{(1)}$ is a regular $V$-adapted process (by Theorem 6.3).

Now, suppose that $I^{(k)}_{V}$, $V\in {\cal P}_{0}$ ,
are well-defined regular $V$-adapted 
processes for which the estimate 
(\ref{8.7}) holds, where $1\leq k \leq n-1$. Then
the integrands 
$X_{C,D}^{\eta}I^{(n-1)}_{E}$
are locally square integrable w.r.t. $dA^{\eta}$ since
$X_{C,D}^{\eta}\in {\cal B}_{{\rm loc}}(D,E)$ and $I_{E}^{(n-1)}$ is regular.
Now, using Proposition 8.1 and condition (\ref{8.5}), we obtain
$$
\sum_{\eta\in {\cal T}}
\int_{0}^{t}
\|
B_{C,D}^{\eta}(s)I^{(n-1)}_{E}(s)
x
\|^{2}d\nu_{u}(s)
$$
$$
\leq \sum_{\eta \in {\cal T}}
\{
\sup_{0\leq r\leq t}
\|B_{C,D}^{\eta}(r)\|^{2}
\}
\{
\max_{0\leq s \leq t}
\|I_{E}^{(n-1)}(s)x\|^{2}
\}
\nu_{u}(t)<\infty
$$
for all $C,D,E\in {\cal P}_{0}$, $t\geq 0$, $w\in {\cal D}_{0}$ and 
$u\in {\cal H}_{0}$. This
implies that $I_{V}^{(n)}$ is a well-defined,
regular $V$-adapted process since each term
in the sum on the RHS of equation (\ref{8.6}) is a regular $V$-adapted process. 
Moreover, using Proposition 8.1 and
then the inductive assumption, we arrive at
$$
\|(I_{V}^{(n)}(t)-I_{V}^{(n-1)}(t))x\|^{2}
$$
\begin{eqnarray*}
&\leq &
|{\cal P}_{0}|^{6}
\max_{C\cap D\cap E= V}
\left\{
\|
\sum_{\eta\in {\cal T}}
\int_{0}^{t}X_{C,D}^{\eta}(I_{E}^{(n-1)}-I_{E}^{(n-2)})
\# dA^{\eta}x 
\|^{2}
\right\}\\
&\leq &
2|{\cal P}_{0}|^{6}
e^{\nu_{u}(T)}
\max_{C\cap D \cap E= V}
\left\{
\sum_{\eta \in {\cal T}}
\sup_{0\leq r \leq t}
\|B_{C,D}^{\eta}(r)\|^{2}\right.\\
&\times &
\left.
\int_{0}^{t}
\|(I_{E}^{(n-1)}(s)-I_{E}^{(n-2)}(s))x\|^{2}d\nu_{u}(s)
\right\}\\
&\leq &
2^{n}|{\cal P}_{0}|^{6n}e^{n\nu_{u}(T)}k_{T}^{n}l_{0}
\|x\|^{2}
\frac{(\nu_{u}(t))^{n}}{n!}.
\end{eqnarray*}
This shows that the estimate (\ref{8.7}) holds for all natural $n$. 
Therefore, the strong limit
$I_{V}(t)x = s-\lim_{n\rightarrow \infty}I_{V}^{(n)}$
exists for all $V\in {\cal P}_{0}$, $t\geq 0$ and $x\in {\cal M}_{0}$
and defines a regular $V$-adapted process which satisfies equation
(\ref{8.4}).

Suppose there exist two solutions $(I_{V})_{V\in {\cal P}_{0}}$,
$(I_{V}')_{V\in {\cal P}_{0}}$ of equation (\ref{8.4}), where
$I_{V}, I_{V}'$ are regular $V$-adapted processes for all $V\in {\cal P}_{0}$.
Then, setting $I_{V}-I_{V}'=Z_{V}$
for all $V\in {\cal P}_{0}$, we have
$$
Z_{V}(t)x
=\sum_{C\cap D \cap E = V}
\sum_{\eta \in {\cal T}}
\int_{0}^{t}
X_{C,D}^{\eta}Z_{E}\#
dA^{\eta} x
$$
for all $x\in {\cal M}_{0}$, $V\in {\cal P}_{0}$ and 
$0 \leq t\leq T$. This gives (by Proposition 8.1)
\begin{eqnarray*}
\max_{V}\|Z_{V}(t)x\|^{2}
&\leq &
|{\cal P}_{0}|^{6}
\max_{V}\max_{C\cap D\cap E=V}
\sum_{\eta \in {\cal T}}
\int_{0}^{t}
\|B_{C,D}^{\eta}(s)Z_{E}(s)x\|^{2}
d\nu_{u}(s)\\
&\leq &
|{\cal P}_{0}|^{6}
\max_{C,D}
\sum_{\eta \in {\cal T}}
\sup_{0\leq r \leq T}\|B_{C,D}^{\eta}(r)\|^{2}
\int_{0}^{t}\max_{E}\|Z_{E}(s)x\|^{2}d\nu_{u}(s)
\end{eqnarray*}
and therefore, by Gronwall's inequality (see [P]), we obtain
$$
\max_{V}\|Z_{V}(t)x\|^{2}\leq 0
$$
for $0\leq t \leq T$,
which gives $Z_{V}(t)=0$ on $\widetilde{\cal E}_{0}$
for $0\leq t \leq T$ and $V\in {\cal P}_{0}$. 
Since $T$ is arbitrary, this implies that
$Z_{V}(t)=0$ for all $t\geq 0$.
\hfill $\Box$\\
\indent{\par}
Using Theorem 8.2, which deals with {\it systems} of
equations, we can establish existence and uniqueness of the solution of 
the stochastic differential {\it equation}
\begin{eqnarray}
\label{8.8}
dI
&=& 
\sum_{\eta \in {\cal T}}
X^{\eta}
I\# dA^{\eta}\\
\label{8.9}
I(0)
&=&
I^{(0)}
\end{eqnarray}
in the class of ${\cal P}_{0}$--adapted processes, 
where $X^{\eta}$ is a suitable
$({\cal P}_{0}, {\cal P}_{0})$-- adapted biprocess
for all $\eta \in {\cal T}$.
We just need to impose a condition
on ${\cal P}_{0}$ which will ensure that processes 
with different types of adaptedness
are linearly independent. \\
\indent{\par}
{\sc Lemma 8.3.}
{\it Let ${\cal P}_{0}=\{V_{1}, \ldots , V_{n}\}$,
be a finite subset of ${\cal P}({\bf N})$
such that}
\begin{equation}
\label{8.10}
V_{i}\setminus (V_{1}\cup \ldots \cup V_{i-1})\neq \emptyset
\end{equation}
{\it for $i=2, \ldots , n$, 
and let $Y_{V}\in {\cal A}(V)$, where $V\in {\cal P}_{0}$.
If}
$$
\sum_{V\in {\cal P}_{0}}Y_{V}(t)x =0
$$
{\it for all $t\geq 0$ and $x\in {\cal M}_{0}$, then 
$Y_{V}(t)x=0$ for all $t\geq 0$,
$x\in {\cal M}_{0}$ and $V\in {\cal P}_{0}$.}\\[5pt]
{\it Proof.}
Let $x=w\varepsilon (u)$. In view of the condition (\ref{8.10}),
$W:=V_{n}$ contains elements which are not in $V_{1}, \ldots , V_{n-1}$.
Let us choose $u$ in such a way that $u^{(W)}_{[t}\neq 0$ 
and $u^{(V)}_{[t}=0$ for all $V\in \{V_{1}, \ldots , V_{n-1}\}$.
We have
$$
\sum_{V\in {\cal P}_{0}}Y_{V}(t)x
=\widetilde{Y}_{W}(t)w\varepsilon (u_{t)})\otimes \varepsilon (u_{[t}^{(W)})
+\sum_{V\neq W}\widetilde{Y}_{V}(t)w\varepsilon (u_{t)}) =0
$$
for all such $u$, $t\geq 0$ and $w\in h_{0}$.
However, by assumption, we also must have
$$
\widetilde{Y}_{W}(t)w\varepsilon (u_{t)})+\sum_{V\neq W}
\widetilde{Y}_{V}(t)w\varepsilon (u_{t)})=0
$$
for all such $u$, $t\geq 0$ and $w\in h_{0}$.
These two facts imply that 
$$
\widetilde{Y}_{W}(t)w\varepsilon (u_{t)})
\otimes \varepsilon_{0} (u_{[t}^{(W)})=0
$$
where 
$$
\varepsilon_{0}(z)=\varepsilon (z) - \Omega
$$
for $z\in {\cal H}_{0}$.
Since $\varepsilon_{0}(u_{[t}^{(W)})\neq 0$, we must have
$$
\widetilde{Y}_{W}(t)w\varepsilon (u_{t)})=0.
$$
Since $u_{t)}$ and $w$ were arbitrary, 
it is now enough to use $W$-adaptedness to get 
$$
Y_{W}(t)w\varepsilon (u) = \widetilde{Y}_{W}(t)w\varepsilon (u_{t)})
\otimes \varepsilon (u_{[t}^{(W)})=0
$$
for all $t\geq 0$, $w\in h_{0}$ and $u\in {\cal H}_{0}$. 
Hence $Y_{W}=0$ on $\widetilde{\cal E}_{0}$.
We can continue this way for other sets in ${\cal P}_{0}$ 
to show that $Y_{V}=0$ for all $V\in {\cal P}_{0}$.
\hfill $\Box$\\
\indent{\par}
In view of Lemma 8.3, the problem of
existence and uniqueness of solutions of 
(\ref{8.8})-(\ref{8.9}) in the 
class of ${\cal P}_{0}$-adapted processes is equivalent to that
of existence and uniqueness of solutions of 
the system of equations (\ref{8.1})-(\ref{8.2}).
The most natural examples of stochastic differential equations 
which ``mix different types of adaptedness'' appear in $m$-free calculi
and will be presented below. Before we do that, let us establish
another result which will be needed in Section 9.

Namely, it is desirable to establish existence and uniqueness of 
solutions for a more general class of equations than 
those given by (\ref{8.8})-(\ref{8.9}). 
In particular, we would like to get uniqueness of solutions of
the equation
\begin{eqnarray}
\label{8.11}
dI &=& \sum_{\eta \in {\cal T}}\{X^{\eta}IM^{\eta}\# dA^{\eta}
+ N^{\eta}IY^{\eta}\# dA^{\eta}\}\\
\label{8.12}
I(0)&=& I^{(0)}
\end{eqnarray}
where 
$$
X^{\eta}=\sum_{D,E\in {\cal P}_{0}}
F^{\eta}_{D}\otimes G^{\eta}_{E},\;\;
Y^{\eta}=\sum_{D,E\in {\cal P}_{0}}
H^{\eta}_{D}\otimes K^{\eta}_{E}
$$
and
$$
M^{\eta}=\sum_{D\in {\cal P}_{0}} M^{\eta}_{D},\;\; 
N^{\eta}=\sum_{D\in {\cal P}_{0}} N^{\eta}_{D}
$$
under suitable assumptions on $X^{\eta},Y^{\eta},M^{\eta},N^{\eta}$.
For our purposes, it will be enough to 
establish conditions under which
the solution of equations (\ref{8.11})-(\ref{8.12})
is unique, if it exists. Note that 
they are not the most general (one can extend 
this result in the spirit of [Ac-Fa-Qu])
and that they also guarantee existence of a solution.

We will say that a $V$-adapted process $F_{V}$
is a $P^{(V)}$--{\it ampliation} if 
$$
F_{V}(t)=\bar{F}_{V}(t)\otimes P^{(V)}
$$
for all $t\geq 0$, according to the decomposition $\widetilde{\Gamma}
=h_{0}\otimes \Gamma$.\\
\indent{\par}
{\sc Theorem 8.4.}
{\it If
$F^{\eta}_{V}$, $H^{\eta}_{V}$, 
$G^{\eta}_{V}$,
$H^{\eta}_{V}$,
$K^{\eta}_{V}$,
$M^{\eta}_{V}$,
$N^{\eta}_{V}$ are $P^{(V)}$-ampliations and}
$$
\bar{F}^{\eta}_{V},\; \bar{G}^{\eta}_{V},\; \bar{H}^{\eta}_{V},\;
\bar{K}^{\eta}_{V},\; \bar{M}^{\eta}_{V}\; \bar{N}^{\eta}_{V} 
\in {\cal B}_{{\rm loc}}(h_{0})
$$
{\it for all $V\in {\cal P}_{0}$ and $\eta\in {\cal T}$, 
then the solution of equations 
(\ref{8.11})-(\ref{8.12}) is unique if it exists.} 
\\[5pt]
{\it Proof.}
Denote by $Z$ the difference of two ${\cal P}_{0}$-adapted
solutions of the above equation. Let
$$
m(t)=\max_{V\in {\cal P}_{0}}\max_{W\in {\cal P}_{0}}\sup_{\|w\|\leq 1}
\|Z_{V}(t)w\varepsilon (u^{(W)})\|
$$
for all $t\geq 0$.
Then, as in the uniqueness proof of Theorem 8.2, 
for $0\leq t \leq T$ and given $u\in {\cal H}_{0}$, there exists
a non-negative constant $c(T,u)$ and a measure $\tau$ which is absolutely 
continuous w.r.t. the Lebesgue measure, such that
$$
m(t)\leq c(T,u)\int_{0}^{t}m(s)d\tau_{s}
$$
for all $0\leq t \leq T$. By Gronwall's inequality, we obtain
$m(t)=0$ on $[0,T]$ and this gives uniqueness of solutions.
This completes the proof of the theorem.
\hfill $\Box$\\[5pt]
\begin{center}
{\it $m$-free stochastic differential equations}
\end{center}
Assume that $m$ is finite and consider stochastic differential equations
of the form 
\begin{eqnarray}
\nonumber
dI
&=&
F_{1}\otimes G_{1}I\# dl^{(m)}
+
F_{2}\otimes G_{2}I\# dl^{(m)*}\\
\label{8.13}
&+&
F_{3}\otimes G_{3}I\# dl^{(m)\circ}
+
F_{4}\otimes G_{4}I\# dl^{(m)\cdot}
\\
\label{8.14}
I(0)&=&I^{(0)}
\end{eqnarray}
where
$$
I^{(0)}=\sum_{V\in {\cal P}_{0}^{(m)}}I_{V}^{(0)}
$$
with $F_{i}\otimes G_{i}\in {\cal B}_{{\rm loc}}({\cal P}_{0},{\cal P}_{0})$
for $1\leq i \leq 4$, where ${\cal P}_{0}^{(m)}$ is given by Example 2
in Section 9 and ${\cal B}_{{\rm loc}}({\cal P}_{0},{\cal P}_{0})$ denotes
the space of locally bounded $({\cal P}_{0},{\cal P}_{0})$-- adapted
biprocesses.

In view of relations (\ref{2.6})-(\ref{2.9}), we can express extended 
$m$-free differentials in terms of filtered differentials and this leads to
\begin{eqnarray}
\nonumber
dI
&=&
\sum_{k=1}^{m}
F_{1}P^{[k-1]}\otimes G_{1}I \# dA^{(k)}
+ 
\sum_{k=1}^{m}
F_{2}\otimes P^{[k-1]}G_{2}I\# dA^{(k)*}\\
\label{8.15}
&+& 
\sum_{k=1}^{m}
F_{3}P^{[k]}\otimes G_{3}I \# dA^{(k)\circ}
+ 
F_{4}P^{(m)}\otimes G_{4}I \# dA^{(0)}
\\[8pt]
\label{8.16}
I(0)&=&I^{(0)}
\end{eqnarray}
Note that equations (\ref{8.15})-(\ref{8.16}) are a
special case of equations (\ref{8.8})-(\ref{8.9}).
From Theorem 8.2 we obtain existence and uniqueness of solutions 
$I_{(m)}$ of equations (\ref{8.13})-(\ref{8.14}) 
for each finite $m\in {\bf N}$.

The case $m=\infty$ is obtained by taking strong limits as we show below.
For that purpose, we need to make stronger assumptions on
the integrated biprocesses $F_{i}\otimes G_{i}$, $1\leq i \leq 4$.\\
\indent{\par}
{\sc Theorem 8.6.}
{\it Suppose that $F_{i}\otimes G_{i}\in {\cal B}_{{\rm loc}}({\cal P}_{0},
{\cal P}_{0})$, $i=1, \ldots , 4$ and that there exists $p\in {\bf N}$
such that the ranges of $F_{i}(t)$ are contained in  
$\widetilde{\Gamma}({\cal H}^{(p)})$ for all $t\geq 0$ and $i=1, \ldots , 4$. 
Then}
\begin{equation}
\label{8.17}
I(t)x=:s-\lim_{m\rightarrow \infty}I_{(m)}(t)x
\end{equation}
{\it exists for all $t\geq 0$ and $x\in {\cal M}_{0}$, and
satisfies equations (\ref{8.13})-(\ref{8.14}) for $m=\infty$.
If $I_{(m)}(t)$ is an isometry for all $t\geq 0$, then $I(t)$
is an isometry for all $t\geq 0$.}\\[5pt]
{\it Proof.}
From the assumption on the ranges of $F_{i}(t)$ and the iteration 
of solutions as in the proof of Theorem 8.2 we obtain
$I_{(m)}^{(n)}(s)x\in \widetilde{\Gamma}({\cal H}^{(p)})$ and thus
$I_{(m)}(s)x\in \widetilde{\Gamma}({\cal H}^{(p)})$
for all $s\geq 0$, $x\in {\cal M}_{0}$ and $m\in {\bf N}$.
Therefore, by Definitions 3.1 and 3.2.,
there exists $q\in {\bf N}$ such that
$$
G_{\sigma}(s)I_{(m)}(s)x, G_{\sigma}(s)I_{(n)}(s) x\in \widetilde{\Gamma}({\cal H}^{(q)})
$$
for all $\sigma =1,2,3,4$, $s\geq 0$ and $n,m\in {\bf N}$. 
From this and relations
(\ref{2.6})-(\ref{2.9}) we infer that there exists $m\in {\bf N}$ such 
that for all $n\geq m$ we have
$$
\int_{0}^{t}F_{\sigma}\otimes G_{\sigma}I_{(n)}x\#dl^{(n)\sigma}=
\int_{0}^{t}F_{\sigma}\otimes G_{\sigma}I_{(n)}x\#dl^{(m)\sigma}
$$
for all $x\in {\cal M}_{0}$, where $\sigma=1,2,3,4$ correspond to
annhilation, creation, number and time processes, respectively.
This leads to
$$
(I_{(m)}(t)-I_{(n)}(t) )x=\sum_{\sigma=1}^{4}
\int_{0}^{t}F_{\sigma}\otimes G_{\sigma}(I_{(m)}-
I_{(n)})x\#dl^{(m)\sigma}.
$$
for all $t\geq 0$ and $n\geq m$.
Decomposing all processes into their filtered components and using
similar arguments as in the proof of Theorem 8.2, we conclude that
$I_{(m)}(t)x=I_{(n)}(t)x$ for $n\geq m$. Therefore, 
$s-\lim_{m\rightarrow \infty}I_{(m)}(t)x=I(t)x$ exists and 
is the unique solution of equation (\ref{8.13})-(\ref{8.14})
for $m=\infty$. The last statement is obvious.
\hfill $\Box$
\\[10pt]
\myownsection
\begin{center}
{\sc 9. Unitary evolutions}
\end{center}
In this section we establish necessary and also sufficient conditions
under which the stochastic differential equation
\begin{eqnarray}
\label{9.1}
dU
&=& 
\sum_{\eta \in {\cal T}_{0}}
X^{\eta}
U\# dA^{\eta}\\
\label{9.2}
U(0)
&=&
1
\end{eqnarray}
has a unique unitary solution $U=(U(t))_{t\geq 0}$, i.e.
$U(t)$ is unitary for all $t\geq 0$, 
where ${\cal P}_{0}$ is a finite subset
of the power set ${\cal P}({\bf N})$, ${\cal T}_{0}$ is a finite subset
of ${\cal T}$ and $X^{\eta}$ are suitable 
$({\cal P}_{0},{\cal P}_{0})$- adapted
biprocesses for all $\eta \in {\cal T}_{0}$.

Throughout this section we assume that ${\cal D}_{0}=h_{0}$ and
we use the notation
$$
B^{\eta}(t)=
\sum_{(D,E)\in {\cal P}_{0}}
\1_{D,E}^{\eta}B_{D,E}^{\eta}(t) 
$$
where $\1_{D,E}^{\eta}$ is given by (4.2) and
$B_{D,E}^{\eta}\models X_{D,E}^{\eta}$ for 
$\eta \in {\cal T}_{0}$ and $D,E\in {\cal P}_{0}$.
We begin with our second linear independence lemma.\\
\indent{\par}
{\sc Lemma 9.1.}
{\it Let $X^{\eta}\in {\cal C}({\cal P}_{0}, {\cal P}_{0})$, where
$\eta \in {\cal T}_{0}$, and 
${\cal T}_{0}$, ${\cal P}_{0}$ are finite subsets of ${\cal T}$ and  
${\cal P}({\bf N})$, respectively, and ${\cal P}_{0}$
is closed under intersections and satisfies the condition (\ref{8.10}).
If\\[3pt]
(i) the map $u\rightarrow  B^{\eta}_{D,E}(t)w\varepsilon (u)$
is strongly continuous for each $t\geq 0$, $w\in h_{0}$,
$u\in {\cal H}_{0}$,
$\eta\in {\cal T}_{0}$,
and $D,E\in {\cal P}_{0}$,\\[5pt]
(ii) $\sum_{\eta \in {\cal T}_{0}}
\int_{0}^{t}X^{\eta}\# dA^{\eta}x =0$
for all $t\geq 0$, $x\in {\cal M}_{0}$, \\[3pt]
then
\begin{equation}
\label{9.3}
B^{\eta}(t)x =0
\end{equation}
for all $\eta \in {\cal T}_{0}$, $x\in {\cal M}_{0}$,
and $t\geq 0$.}\\[5pt]
{\it Proof.}
The proof is a straightforward modification of that for adapted processes
[Par] and therefore will be omitted.\hfill $\Box$\\
\indent{\par}
Let us now address the question of unitarity of the solution
of (\ref{9.1})-(\ref{9.2}).
In other words, we are looking for necessary and sufficient
conditions under which $U(t)U^{*}(t)=U^{*}(t)U(t)=1$ 
for all $t\geq 0$. If ${\cal P}_{0}$ is closed under intersections
and the condition (\ref{8.10}) is satisfied, then,
for the solution to be unitary it is necessary that
${\bf N}\in {\cal P}_{0}$ since we then must have
$$
\sum_{C,D\in {\cal P}_{0}}U_{C}^{*}(t)U_{D}(t)=
\sum_{V\in {\cal P}_{0}}\sum_{C\cap D=V}U_{C}^{*}(t)U_{D}(t)=1
$$
$$
\sum_{C,D\in {\cal P}_{0}}U_{C}(t)U_{D}^{*}(t)=
\sum_{V\in {\cal P}_{0}}\sum_{C\cap D=V}U_{C}(t)U_{D}^{*}(t)=1
$$
for all $t\geq 0$. This, by Lemma 8.3, 
implies that the sum over $V$ must include
$V={\bf N}$ and that
$$
U_{{\bf N}}^{*}(t)U_{{\bf N}}(t)=U_{{\bf N}}(t)U_{{\bf N}}^{*}(t)=1
$$
for all $t\geq 0$, whereas
$$
\sum_{C\cap D=V}
U_{C}^{*}(t)U_{D}(t)=\sum_{C\cap D=V}U_{C}(t) U_{D}^{*}(t)=0
$$
for all $t\geq 0$ and any $V\in {\cal P}_{0}$ such that $V\neq {\bf N}$.
This means that in order to study unitarity on $\Gamma({\cal H})$
one needs to include adapted biprocesses in the
filtered adapted biprocesses.
For that reason we will assume that
${\cal P}_{0}$ contains ${\bf N}$.
However, in order to establish unitarity conditions, stronger
conditions on ${\cal P}_{0}$ are needed.
This motivates the following definition. \\
\indent{\par}
{\sc Definition 9.2.}
We will say that a collection 
${\cal P}_{0}\subset {\cal P}({\bf N})$
is {\it admissible} if 
$$
{\cal P}_{0}=\{V_{i},  1\leq i \leq p \},\;\; 1\leq p \leq \infty
$$
where $V_{i}$ is a proper subset of $V_{i+1}$ for all $1\leq i\leq p$ and 
$V_{p}={\bf N}$.\\
\indent{\par}
{\it Example 1.}
Let ${\cal P}_{0}=\{{\bf N}\}$, i.e. ${\cal P}_{0}$ consists 
of one filter which corresponds to adapted biprocesses,
i.e. {\it boson calculus} (either with a finite or infinite number
of degrees of freedom). Then ${\cal P}_{0}$ is admissible.\\
\indent{\par}
{\it Example 2.}
Let ${\cal P}_{0}^{(m)}=\{V(k), {\bf N}: 1\leq k \leq m +1 \}$, where
$m\in {\bf N}$, with $V(k)=\{1, \ldots , k-1\}$. 
Then ${\cal P}_{0}^{(m)}$ is finite and admissible for each $m\in {\bf N}$.
These collections of filters appear in $m$-{\it free calculi}
for finite $m$. In turn ${\cal P}_{0}^{(\infty)}=
\{V(k), {\bf N}: 1\leq k < \infty\}$ is an admissible 
collection of filters corresponding to {\it free calculus}.\\
\indent{\par}
{\sc Lemma 9.3.}
{\it Let ${\cal P}_{0}\subset {\cal P}({\bf N})$ be admissible
and finite and suppose that 
$X^{\eta}\in {\cal B}_{{\rm loc}}({\cal P}_{0},{\cal P}_{0})$
$\cap {\cal C}({\cal P}_{0},{\cal P}_{0})$ and are non-zero
if and only if $\eta\in {\cal T}_{0}$, where
${\cal T}_{0}$ is a finite subset of ${\cal T}$ and that
the continuity condition (i) of Lemma 9.1 is satisfied. 
Then, for the unique solution of equations (\ref{9.1})-(\ref{9.2} )
to be an isometry it is necessary that} 
$$
(B^{\eta\dagger}(t))^{*}+B^{\eta}(t) 
+
\sum_{\bl \eta_{1}, \eta_{2} \br =\eta}
(B^{\eta_{1}\dagger}(t))^{*}B^{\eta_{2}}(t)=0
$$
{\it for all $t\geq 0$ and $\eta \in {\cal T}$, where}
\begin{equation}
\label{9.4}
\sum_{\bl \eta_{1}, \eta_{2} \br =\eta}
:=
\sum_{
\stackrel{\eta_{1}, \eta_{2}}
{\scriptscriptstyle \bl A^{\eta_{1}} , A^{\eta_{2}} \br = A^{\eta}}
}
\end{equation}
{\it Proof.}
The proof is based on the filtered It$\hat{{\rm o}}$ formula. 
We will use elementary biprocesses, i.e. 
$X^{\eta}_{D,E}=F_{D}^{\eta}\otimes G_{E}^{\eta}$ for all
$\eta \in {\cal T}_{0}$ and $D,E\in {\cal P}_{0}$.
In all summations it is implicitly assumed that
$\eta \in {\cal T}_{0}$ and $D,E\in {\cal P}_{0}$
and only additional conditions on the summation indices
are shown. We have
\begin{eqnarray*}
dU
&=&
\sum_{\eta}\sum_{D,E}
F_{D}^{\eta}\otimes G_{E}^{\eta}U\# dA^{\eta}\\
U(0)
&=&
1
\end{eqnarray*}
and since the number of terms in the sum is finite,
it is clear that $U^{*}=(U^{*}(t))_{t\geq 0}$ satisfies
\begin{eqnarray*}
dU^{*}(t)
&=&
\sum_{\eta}\sum_{D,E}
U^{*}(G_{E}^{\eta})^{*} 
\otimes 
(F_{D}^{\eta})^{*} 
\# 
dA^{\eta\dagger}\\
U^{*}(0)
&=&
1
\end{eqnarray*}
Applying the filtered It$\hat{{\rm o}}$ formula to the isometry condition
$U^{*}(t)U(t)=1$, we obtain (skipping $t$ to save some space and
choosing one way of writing the It$\hat{{\rm o}}$ correction of Theorem 6.1)
\begin{eqnarray*}
0 &=&
dU^{*}U+U^{*}dU + dU^{*}dU\\
&=&
\sum_{\eta}\sum_{D,E}U^{*}
(G_{E}^{\eta})^{*}
\otimes
(F_{D}^{\eta})^{*}U 
\# dA^{\eta \dagger}
+
\sum_{\eta}\sum_{D,E}U^{*}F_{D}^{\eta}\otimes
G_{E}^{\eta}U \# dA^{\eta}\\
&+&
\sum_{\eta_{1},\eta_{2}} 
\sum_{D_{1},E_{1},D_{2},E_{2}}
U^{*}
(G_{E_{1}}^{\eta_{1}\dagger})^{*} 
\otimes
\rho_{\eta_{1},\eta_{2}}
\left[
(F_{D_{1}}^{\eta_{1}\dagger})^{*} 
F_{D_{2}}^{\eta_{2}}
\right]
G_{E_{2}}^{\eta_{2}}U \# 
d\bl A^{\eta_{1}}, A^{\eta_{2}} \br
\end{eqnarray*}
Note that $U$ and $U^{*}$ are ${\cal P}_{0}$-adapted, i.e.
in general, they contain mixed types of adaptedness. Using
equation (\ref{9.3}), we obtain
\begin{equation}
\label{9.5}
\sum_{(W,Z)\in {\cal P}_{0}}
\1_{W,Z}^{\eta}U_{W}^{*}
\{
(B^{\eta \dagger})^{*}
+
B^{\eta}
+
\sum_{\bl \eta_{1}, \eta_{2} \br =\eta}
(B^{\eta_{1}\dagger})^{*}
B^{\eta_{2}}
\}
U_{Z}=0
\end{equation}
for each $\eta \in {\cal T}_{0}$. 

Since ${\cal P}_{0}$ is increasing, it is easy to check that we have
$$
\sum_{
\stackrel{V\in {\cal P}_{0}}
{\scriptscriptstyle k\in V}
}
U_{V}(t)
\sum_{
\stackrel{W\in {\cal P}_{0}}
{\scriptscriptstyle k\in Z}
}
U_{W}^{*}(t)
=
U_{{\bf N}}(t)U_{{\bf N}}^{*}(t)
=1
$$
and therefore, by 
mulitplying equation (\ref{9.5}) by $U$ or 
$
\sum_{
\stackrel{V\in {\cal P}_{0}}
{\scriptscriptstyle k\in V}
}
U_{V}(t)
$
from the left and by $U^{*}$ or
$
\sum_{
\stackrel{V\in {\cal P}_{0}}
{\scriptscriptstyle k\in V}
}
U_{V}^{*}(t)
$
from the right (that depends on $\eta$), we arrive at
$$
(B^{\eta\dagger}(t))^{*}+B^{\eta}(t) 
+
\sum_{\bl \eta_{1}, \eta_{2} \br =\eta}
(B^{\eta_{1}\dagger}(t))^{*}B^{\eta_{2}}(t)=0
$$
which ends the proof.
\hfill
$\Box$\\
\indent{\par}
{\sc Lemma 9.4.}
{\it Under the assumptions of Lemma 9.3, for the unique solution of
equations (\ref{9.1})-(\ref{9.2}) to be a co-isometry it is necessary that}
$$
(B^{\eta\dagger}(t))^{*}+B^{\eta}(t) 
+
\sum_{\bl \eta_{1}, \eta_{2} \br =\eta}
B^{\eta_{1}}(t)(B^{\eta_{2}\dagger}(t))^{*}=0
$$
{\it for all $t\geq 0$ and $\eta \in {\cal T}_{0}$.}\\[5pt] 
{\it Proof.}
The proof is similar to that of Lemma 9.3
and is based on differentiating the co-isometry
condition $U(t)U(t)^{*}=1$ and then
using the filtered It$\hat{{\rm o}}$ formula. \hfill $\Box$\\
\indent{\par}
{\sc Theorem 9.5.}
{\it Under the assumptions of Lemma 9.3,
for the unique solution of equations (\ref{9.1})-(\ref{9.2}) to be unitary, it
is necessary that for all $t\geq 0$}
$$
\begin{array}{ll}
(i) & B^{(k)\circ}(t)+1\;\;{\rm is}\;\; {\rm unitary}
\;\; \forall k\in {\bf N}\\[5pt]
(ii)& (B^{(k)*}(t))^{*} + 
B^{(k)}(t) + (B^{(k)*}(t))^{*}B^{(k)\circ}(t) =0\\[5pt]
(iii)& B^{(0)}(t) + (B^{(0)}(t))^{*} +
\sum_{k\geq 1}(B^{(k)*}(t))^{*} B^{(k)*}(t)=0
\end{array}
$$
{\it Proof.}
It is a straightforward consequence of Lemmas 9.3-9.4.
\hfill $\Box$\\
\indent{\par}
{\it Remark.}
Although the unitarity conditions of Theorem 9.5 have the same form as
in boson calculus, it is important that 
$B^{\eta}$ are, in general, ${\cal P}_{0}$-adapted and
not ${\bf N}$-adapted processes.
Note also that, in general, certain summands in (\ref{9.4}) 
are equal to zero.
This tells us which components of $({\cal P}_{0},{\cal P}_{0})$-adapted
biprocesses may give a non-zero contribution
to the differential equation (only these enter the unitarity
conditions).\\
\indent{\par}
{\sc Theorem 9.6.}
{\it Suppose ${\cal P}_{0}$ is admissible and finite and 
$X_{D,E}^{\eta}=F_{D}^{\eta}\otimes G_{E}^{\eta}$, with
$F^{\eta}_{D}(t)=\bar{F}^{\eta}_{D}(t)\otimes P^{(D)}$,
$G^{\eta}_{E}(t)=\bar{G}^{\eta}_{E}(t)\otimes P^{(E)}$
according to the decomposition $\widetilde{\Gamma}=h_{0}\otimes
\Gamma$, where
$\bar{F}^{\eta}_{D}, \bar{G}^{\eta}_{E}\in 
{\cal B}_{{\rm loc}}(h_{0})$ for
all $\eta \in {\cal T}_{0}$ and $D,E\in {\cal P}_{0}$.
Then the conditions (i)-(iii) of Theorem 9.5 are sufficient for the unique
solution of equations (\ref{9.1})-(\ref{9.2}) to be unitary.}\\[5pt]
{\it Proof.}
Let 
$$
U(t)=1+\sum_{\eta}\int_{0}^{t}X^{\eta}
U\# dA^{\eta}
$$
for all $t\geq 0$. We will first show that if conditions (i)-(iii)
of Theorem 9.5 are satisfied, 
then $U$ is an operator-valued isometric process. Denote
$x=w\varepsilon (u) , z= y\varepsilon (v)$, where $w,y\in h_{0}$
and $u,v\in {\cal H}_{0}$ and let
$$
I_{Z}(t)=\sum_{\eta}\sum_{D\cap E\cap V=Z}
\int_{0}^{t}F_{D}^{\eta}\otimes G_{E}^{\eta}
U_{V}\# dA^{\eta}
$$
for $Z\in {\cal P}_{0}$ and $t\geq 0$.
We have 
\begin{eqnarray*}
\langle U(t)x, U(t)y \rangle
-
\langle x , y \rangle
&=&
\sum_{\eta}\sum_{W}\int_{0}^{t}
\langle U_{W}x , 
(B^{\eta})^{*}y \rangle d\mu_{W, {\bf N}}^{\eta \dagger}\\
&+&
\sum_{\eta}\sum_{W}\int_{0}^{t}
\langle x, 
B^{\eta}U_{W}y 
\rangle
d\mu_{{\bf N},W}^{\eta}\\
&+&
\sum_{\eta_{1}}\sum_{W_{1},Z}
\int_{0}^{t}
\langle
U_{W_{1}}x,
(B^{\eta_{1}})^{*}
I_{Z}y
\rangle
d\mu_{W_{1} , Z}^{\eta_{1}\dagger}\\
&+&
\sum_{\eta_{2}}
\sum_{W_{2},Z}
\int_{0}^{t}
\langle
I_{Z}x,
B^{\eta_{2}}U_{W_{2}} y
\rangle
d\mu_{Z, W_{2}}^{\eta_{2}}\\
&+&
\sum_{\eta}
\sum_{\bl \eta_{1},\eta_{2} \br =\eta}
\sum_{W_{1},W_{2}}
\int_{0}^{t}
\langle
U_{W_{1}}x,
(B^{\eta_{1}\dagger})^{*}
B^{\eta_{2}}U_{W_{2}}y
\rangle
d\mu_{W_{1}, W_{2}}^{\eta}
\end{eqnarray*}
where we used the filtered It$\hat{{\rm o}}$ formula. Now, using
\begin{eqnarray*}
U_{Z}&=&I_{Z}, \;\;\; Z\neq {\bf N}\\
U_{{\bf N}}&=&1 + I_{{\bf N}},
\end{eqnarray*}
then 
replacing $Z$ in the 3rd and 4th
terms by $W_{2}$ and $W_{1}$, respectively, and
taking into account all cancellations,
we arrive at
\begin{eqnarray*}
&& \sum_{\eta}\sum_{W_{1},W_{2}}
\int_{0}^{t}
\langle
U_{W_{1}}x, (B^{\eta\dagger})^{*} U_{W_{2}} y
\rangle
d\mu_{W_{1},W_{2}}^{\eta}\\
&+&
\sum_{\eta}\sum_{W_{1},W_{2}}
\int_{0}^{t}
\langle
U_{W_{1}}x, B^{\eta} U_{W_{2}} y
\rangle
d\mu_{W_{1},W_{2}}^{\eta}\\
&+&
\sum_{\eta}
\sum_{\bl \eta_{1},\eta_{2} \br =\eta}
\sum_{W_{1},W_{2}}
\int_{0}^{t}
\langle
U_{W_{1}}x, (B^{\eta_{1}\dagger})^{*}B^{\eta_{2}} U_{W_{2}} y
\rangle
d\mu_{W_{1},W_{2}}^{\eta}=0
\end{eqnarray*}
in view of the isometry conditions of Lemma 9.3. 
Therefore, $U$ is an operator-valued isometric process.

Proceeding in a similar way, we can show that if, say 
$$
{\cal P}_{0}=\{V_{1}, \ldots , V_{n}\},\;\;{\rm where}\;\; 
V_{1}\subset V_{2}\subset \ldots \subset V_{n}={\bf N}
$$
then
$$
\sum_{i=k}^{n}U_{V_{i}}(t)
$$
is also an isometry for all $1\leq k \leq n$ and $t\geq 0$. Therefore,
each $U_{V}(t)$, $V\in {\cal P}_{0}$, is an operator-valued process.

Let us show that $U(t)$ is a co-isometry for all $t\geq 0$.
The adjoint process $U^{*}(t)$ obeys the equation
\begin{eqnarray*}
dU^{*}(t)&=&\sum_{\eta}\sum_{D,E}U^{*}(t)
(G_{E}^{\eta}(t))^{*}
\otimes (F_{D}^{\eta}(t))^{*}\# dA^{\eta\dagger}_{t}\\
U^{*}(0)&=&1.
\end{eqnarray*}
Proceeding as in the isometry case, we obtain
\begin{eqnarray*}
\langle
U^{*}(t)x, U^{*}(t)y
\rangle
-
\langle
x , y 
\rangle
&=&
\sum_{\eta}\sum_{W_{1},W_{2}}
\int_{0}^{t}
\langle
x,
U_{W_{1}}U_{W_{2}}^{*}(B^{\eta\dagger})^{*}y
\rangle
d\mu_{W_{1}\cap W_{2}, {\bf N}}^{\eta}\\
&+&
\sum_{\eta}\sum_{W_{1},W_{2}}
\int_{0}^{t}
\langle
x, B^{\eta}
U_{W_{1}}U_{W_{2}}^{*} y
\rangle
d\mu_{{\bf N}, W_{1}\cap W_{2}}^{\eta}\\
&+&
\sum_{\eta}
\sum_{\bl \eta_{1} , \eta_{2} \br =\eta}
\sum_{W_{1},W_{2}}
\int_{0}^{t}
\langle
x,
B^{\eta_{1}}U_{W_{1}}U_{W_{2}}^{*}(B^{\eta_{2}\dagger})^{*}
\rangle
d\mu_{W_{1}\cap W_{2},W_{1}\cap W_{2}}^{\eta}.
\end{eqnarray*}
This equation is equivalent to the stochastic differential equation
\begin{eqnarray*}
d(UU^{*})&=&
\sum_{\eta}\sum_{D,E}
\{
UU^{*}(G_{E}^{\eta\dagger})^{*}
\otimes (F_{D}^{\eta\dagger})^{*}\# dA^{\eta}
+
F_{D}^{\eta}\otimes G_{E}^{\eta}UU^{*}
\#
dA^{\eta}
\}
\\
&+&
\sum_{\eta}\sum_{\bl \eta_{1},\eta_{2} \br =\eta}
\sum_{
\stackrel{D_{1},E_{1}}
{\scriptscriptstyle D_{2},E_{2}}
}
F_{D_{1}}^{\eta_{1}}
\otimes
G_{E_{1}}^{\eta_{1}}UU^{*}
(G_{E_{2}}^{\eta_{2}\dagger})^{*}(F_{D_{2}}^{\eta_{2}\dagger})^{*}
\# dA^{\eta}
\end{eqnarray*}
Now, $U(t)U^{*}(t)=1$ is a solution of this equation if the
co-isometry conditions of Lemma 9.4 hold. This solution is unique in view
of Theorem 8.4, which completes the proof.
\hfill $\Box$
\begin{center}
{\it Unitarity conditions for boson calculus}
\end{center}
\indent{\par} 
Note that the unitarity conditions of Theorem 9.5 have the same form 
as those for boson calculus on multiple symmetric Fock space
[Mo-Si] (cf. [H-P1], see also [P]). To recover them, 
set ${\cal P}_{0}=\{{\bf N}\}$, $B^{(k)*}=L_{k}$,
$B^{(k)\circ}=S_{k}-1$ and $R=H$. \\[5pt]
\begin{center}
{\it Unitarity conditions for $m$-free calculi}
\end{center}
Let us show that Theorem 9.5 covers unitarity conditions
for $m$-free calculi.
Let $m\in {\bf N}$ and let ${\cal P}_{0}={\cal P}_{0}^{(m)}$ 
(see Example 2 in this section).

For a given ${\cal P}_{0}$-adapted process 
$F=\sum_{V\in {\cal P}_{0}}F_{V}$, let 
$$
[F]_{k}=
\sum_{
\stackrel{V\in {\cal P}_{0}
}
{\scriptscriptstyle k\in V}
}
F_{V}
$$
where $k\in {\bf N}$. Using this notation, we can write
\begin{eqnarray*}
B^{(k)}(t)&=& F_{1}(t)P^{[k-1]}[G_{1}(t)]_{k}\\
B^{(k)*}(t)&=& [F_{2}(t)]_{k}P^{[k-1]}G_{2}(t)\\
B^{(k)\circ}(t)&=& [F_{3}(t)]_{k}P^{(k+1)}[G_{3}(t)]_{k}\\
B^{(0)}(t)&=&F_{4}(t)P^{(m)}G_{4}(t).
\end{eqnarray*}
It is important to notice that 
$[P^{[k]}F]_{k}=P^{(k+1)}[F]_{k}$
which gives the third equation above. The other ones just express
relations between two notations.

In particular, when $F_{3}(t)=G_{3}(t)= 0$ and $F_{2}(t)=G_{1}(t)=1$
for all $t\geq 0$, the above conditions can be written
in an equivalent form
\begin{eqnarray}
F_{4}(t)P^{(m)}G_{4}(t)+G_{4}^{*}(t)P^{(m)}F_{4}^{*}(t) + 
G_{2}^{*}(t)G_{2}(t)&=&0\\
F_{1}(t)P^{(m)}+ G_{2}^{*}(t)P^{(m)}&=&0
\end{eqnarray}
which are the $m$-truncated versions of the unitarity conditions
in [K-Sp]. 

From the considerations of Section 8 it follows that if 
$F_{i},G_{i}$, $i=1, \ldots , 4$, 
are linear combinations of $P^{(V)}$-ampliations for 
$V\in {\cal P}_{0}$ which satisfy the assumptions of Theorem 9.6, then
the strong limit of unitary solutions $s$--$\lim_{m\rightarrow \infty}U_{(m)}$
exists and is unitary. Therefore, it is sufficient that
the $m$-truncated unitarity conditions hold for all $m\in {\bf N}$ which 
is equivalent to 
\begin{eqnarray*}
F_{4}(t)G_{4}(t) + G_{4}^{*}(t) F_{4}^{*}(t)+ 
G_{2}^{*}(t) G_{2}(t)&=&0\\
F_{1}(t) + G_{2}^{*}(t)& =&0
\end{eqnarray*}
i.e. the unitarity conditions for the free calculus.
\\[10pt]
\begin{center}
{\sc Acknowledgements}
\end{center}
I would like to thank Professor Cz. Ryll-Nardzewski, Professor V. Belavkin
and Professor M. Lindsay for helpful remarks.
\begin{center}
{\sc References}\\[20pt]
\end{center}
[Ac-Fa-Qu] {\sc L.~Accardi, F.~Fagnola, J.~Quaegebeur}, 
A representation free quantum stochastic calculus, 
{\it J. Funct. Anal.} {\bf 104} (1992), 149-197.\\[3pt]
[At-Lin] {\sc S.~Attal, M.~Lindsay}, Quantum stochastic calculus --
a maximal formulation, Preprint No. 456 (1999), CNRS, Univ.
de Grenoble I.\\[3pt]
[Ap-H] {\sc D.~B.~Applebaum, R.~H.~Hudson}, Fermion It$\hat{{\rm o}}$'s
formula and stochastic evolutions, {\it Commun.~Math.~Phys.}
{\bf 96} (1984), 437-496.\\[3pt]
[B] {\sc A.~Boukas}, Quantum stochastic Analysis: a non-Brownian
case, Ph.D Thesis, Southern Illinois University, U.S.A. 1988.\\[3pt]
[Ba-St-Wi] {\sc C.~Barnett, R.~F.~Streater, I.~F.~Wilde},
Quasi-free quantum stochastic integrals for the CAR and CCR, 
{\it J.~Funct.~Anal.} {\bf 52} (1983), 19-47.\\[3pt]
[Be] {\sc A.~Belton}, Quantum $\Omega$-semimartingales and stochastic 
evolutions, preprint, 2000.\\[3pt]
[Bel] {\sc V.~P.~Belavkin}, Chaotic states and stochastic integration
in quantum systems, {\it Uspekhi Mat. Nauk} {\bf 47} (1992), 47-106.\\[3pt]
[Bi-Sp] {\sc P.~Biane P., R.~Speicher},
Stochastic calculus with respect to free Brownian motion and analysis
on Wigner space, preprint.\\[3pt]
[Fa] {\sc F.~Fagnola}, On quantum stochastic integration with respect
to free noises, {\it in} ``Quantum Probability and Related Topics'',
World Scientific, London, 1991.\\[3pt]
[F-L] {\sc U.~Franz, R.~Lenczewski},
Limit theorems for the hierarchy of freeness, 
{\it Prob. Math. Stat.} {\bf 19} (1999), 23-41.\\[3pt]
[F-L-S] {\sc U.~Franz, R.~Lenczewski, M.~Sch\"{u}rmann}, 
The GNS construction for the hierarchy of freeness, Preprint No. 9/98,
Wroclaw University of Technology, 1998.\\[3pt]
[H-P1] {\sc R.~L.~Hudson, K.~R.~Parthasarathy},
Quantum It$\hat{{\rm o}}$'s formula and stochastic evolutions, 
{\it Commun.~Math.~Phys.} {\bf 93}, 301-323 (1984).\\[3pt]
[H-P2] {\sc R.~L.~Hudson, K.~R.~Parthasarathy}, Unification of
Fermion and Boson stochastic calculus, {\it Commun.~Math.~Phys.}
{\bf 104} (1986), 457-470.\\[3pt]
[K-Sp] {\sc B.~K\"{u}mmerer, R.~Speicher},
Stochastic integration on the Cuntz Algebra ${\cal O}_{\infty}$,
{\it J.~Funct.~Anal.} {\bf 103} (1992), 372-408.\\[3pt]
[L1] {\sc R.~Lenczewski}, Unification of independence in
quantum probability, {\it Inf.~Dim. Anal.~Quant.~Probab. \& Rel.~Top.}
{\bf 1} (1998), 383-405.\\[3pt]
[L2] {\sc R.~Lenczewski}, Filtered random 
variables, preprint, Wroc{\l}aw University of Technology.\\[3pt]
[Li] {\sc M.~Lindsay}, Quantum and non-causal stochastic calculus,
{\it Probab.~Th.~Rel.~Fields} {\bf 97} (1993), 65-80.\\[3pt]
[Ma] {\sc H.~Maassen}, Quantum Markov processes on Fock space described
by integral kernels, {\it in} ``Quantum probability and applications II'',
Eds. L.~Accardi, W.~v.Waldenfels, Lecture Notes in Math., vol. 1136,
Springer, Berlin, 1985.\\[3pt]
[Me] {\sc P.~A.~Meyer}, ``Quantum probability for probabilists'',
Lecture Notes in Math., Vol.~1538, Springer-Verlag, New York/Berlin,
1993.\\[3pt]
[Mo-Si] {\sc A.~Mohari, K.~B.~Sinha}, Quantum stochastic flows 
with infinite degrees of freedom and countable state Markov processes,
{\it Sankhya}, Ser.A, Part I (1990), 43-57.\\[3pt]
[P] {\sc K.~R.~Parthasarathy}, ``An Introduction to Quantum
Stochastic Calculus'', Birkha\"{a}user, Basel, 1992. \\[3pt]
[Par-Sin] {\sc K.~R.~Parthasarathy, K.~B.~Sinha}, Unification of
quantum noise processes in Fock spaces,
Proc. Trento conference on Quantum Probability and Applications
(1989).\\[3pt]
[S1] {\sc M.~Sch\"{u}rmann}, White Noise on Bialgebras, Lecture Notes
in Math. 1544, Springer, Berlin, 1991.\\[3pt] 
[S2] {\sc M.~Sch\"{u}rmann}, Direct sums of tensor products and
non-commutative independence, {\it J. Funct. Anal.} {\bf 133}
(1995), 1-9.\\[3pt]
[Sp1] {\sc R.~Speicher}, On universal products, {\it Fields Institute
Commun.} {\bf 12} (1997), 257-266.\\[3pt]
[Sp2] {\sc R.~Speicher}, Stochastic integration on the full Fock space 
with the help of a kernel calculus, {\it Publ.~R.I.M.S. Kyoto Univ.}
{\bf 27} (1991), 149-184.\\[3pt]
[Vi] {\sc G.~F.~Vincent-Smith}, Classical and quantum dynamical
propagators, Oxford, preprint, 1990.\\[3pt]
[V] {\sc D.~Voiculescu}, Symmetries of some reduced free product 
${\cal C}^{*}$-algebras, {\it in} ``Operator Algebras and their 
Connections with Topology and Ergodic Theory'', Lecture
Notes in Math. 1132, Springer, Berlin, 1985, 556-588.\\[3pt]
\end{document}